\documentclass[10pt]{amsart}

\usepackage{amssymb,amsmath,amsthm,amscd}

\advance\oddsidemargin by -1cm
\advance\evensidemargin by -1cm
\newtheorem{theorem}{Theorem}
\newtheorem*{theorem1}{Theorem 1}
\newtheorem{proposition}{Proposition}
\newtheorem{corollary}{Corollary}
\newtheorem{lemma}{Lemma}
\newtheorem{definition}{Definition}

\theoremstyle{definition}
\newtheorem{remark}{Remark}

\newcommand{\bdm}{\begin{displaymath}}
\newcommand{\edm}{\end{displaymath}}
\newcommand{\bq}{\begin{equation}}
\newcommand{\eq}{\end{equation}}
\newcommand{\bqn}{\begin{equation*}}
\newcommand{\eqn}{\end{equation*}}

\newcommand{\norm}[1]{\left\| #1 \right\|}
\newcommand{\mklm}[1]{\left\{ #1 \right\}}
\newcommand{\eklm}[1]{\left\langle #1 \right\rangle}

\renewcommand{\d}{\,d}
\newcommand{\N}{{\mathbb N}}

\newcommand{\C}{{\mathbb C}}
\newcommand{\R}{{\mathbb R}}

\newcommand{\D}{{\mathcal D}}
\newcommand{\E}{{\mathcal E}}
\newcommand{\F}{{\mathcal F}}

\renewcommand{\H}{{\mathcal H}}
\newcommand{\I}{{\mathcal I}}
\newcommand{\M}{{\mathcal M}}
\newcommand{\T}{{\rm T}}

\renewcommand{\O}{{\mathcal O}}

\newcommand{\1}{{\bf 1}}
\renewcommand{\epsilon}{\varepsilon}
\renewcommand{\phi}{\varphi}
\renewcommand{\rho}{\varrho}

\newcommand{\Cinft}{{\rm C^{\infty}}}
\newcommand{\CT}{{\rm C^{\infty}_c}}

\renewcommand{\L}{{\rm L}}
\newcommand{\Lcal}{{\mathcal L}}
\newcommand{\Ncal}{{\mathcal N}}

\renewcommand{\S}{{\mathcal S}}

\newcommand{\GL}{\mathrm{GL}}

\renewcommand{\O}{{\mathrm O}}

\renewcommand{\det}{\mathrm{det}\,}

\newcommand{\vol}{\text{vol}\,}
\newcommand{\dist}{\text{dist}\,}
\newcommand{\Op}{\mathrm{Op}}

\DeclareMathOperator{\Ker}{Ker}
\DeclareMathOperator{\supp}{supp}

\DeclareMathOperator{\tr}{tr}
\DeclareMathOperator{\gd}{\partial}

\newcommand{\dbar}{{\,\raisebox{-.1ex}{\={}}\!\!\!\!d}}

\begin{document}

\author{Pablo Ramacher}
\title[Reduced Weyl asymptotics for  PDO on bounded domains I]{Reduced Weyl asymptotics for pseudodifferential operators on bounded domains I  \\ The finite group case} 
\address{Pablo Ramacher, Georg-August-Universit\"at G\"ottingen, Institut f\"ur Mathematik, Bunsenstr. 3-5, 37073 G\"ottingen, Germany}
\email{ramacher@uni-math.gwdg.de}
\thanks{The author was supported by the grant RA 1370/1-1 of the German Research Foundation (DFG) during the preparation of this work.}

\begin{abstract} Let $G\subset \O(n)$ be a  group of isometries acting on $n$-dimensional Euclidean space $\R^n$, and ${\bf{X}}$ a bounded  domain in $\R^n$ which is transformed into itself under the action of $G$. Consider a  symmetric, classical pseudodifferential operator $A_0$ in $\L^2(\R^n)$ with $G$-invariant Weyl symbol, and assume that it is semi-bounded from below. We show that the spectrum of  the Friedrichs extension $A$ of the operator $\mathrm{res} \circ A_0 \circ \mathrm{ext}: \CT({\bf{X}}) \rightarrow \L^2({\bf{X}})$ is discrete, and derive asymptotics for the number  $N_\chi(\lambda)$ of eigenvalues of $A$  less or equal $\lambda$ and with eigenfunctions in the $\chi$-isotypic component of $\L^2({\bf{X}})$, giving also an estimate for the remainder term in both cases where $G$ is a finite, or, more generally, a compact group. In particular, we show that the multiplicity of each unitary irreducible representation in $\L^2({\bf{X}})$ is asymptotically proportional to its dimension.
\end{abstract}

\maketitle

\section{Statement of the problem}

Let $G\subset \O(n)$ be a compact group of isometries acting on Euclidean space $\R^n$, and ${\bf{X}}$ a bounded  domain in $\R^n$  which is transformed into itself  under the action of $G$. Consider the regular  representation of $G$
\begin{equation*}
\T(g) \phi(x)=\phi(g^{-1}x)
\end{equation*}
in the Hilbert spaces $\L^2(\R^n)$, and $\L^2({\bf{X}})$, respectively, and endow them with some $G$-invariant scalar product $(\cdot,\cdot)$, so that the representation $T$ becomes unitary. As a consequence of the Peter-Weyl Theorem,  the representation $T$ decomposes into isotypic components according to 
\bqn 
\L^2(\R^n)=\bigoplus _{\chi \in \hat G} \H_\chi, \qquad \L^2({\bf{X}})= \bigoplus _{\chi \in \hat G} \mathrm{res} \, \H_\chi,
\eqn
where $\hat G$ denotes the set of irreducible characters of $G$,  and $\mathrm{res}:\L^2(\R^n) \rightarrow \L^2({\bf{X}})$ is the natural restriction operator. Similarly, $\mathrm{ext}: \CT({\bf{X}}) \rightarrow \L^2(\R^n)$ will denote the natural extension operator. Let $A_0$ be a  symmetric, classical pseudodifferential operator in $\L^2(\R^n)$ of order $2m$ with $G$-invariant Weyl symbol $a$  and principal symbol $a_{2m}$, and assume that $(A_0 \, u, \, u) \geq c \norm{u}^2_m$ for some $c>0$ and  all $u \in \CT({\bf{X}})$, where $\norm{\cdot}_s$ is a norm in the Sobolev space $\mathrm{H}^s(\R^n)$. Consider further  the Friedrichs extension of the  lower semi-bounded operator
\bqn
\mathrm{res} \circ A_0 \circ \mathrm{ext}: \CT({\bf{X}}) \longrightarrow \L^2({\bf{X}}),
\eqn
which is a self-adjoint operator in $\L^2({\bf{X}})$, and denote it by $A$. Finally, let $\gd {\bf{X}}$ be the boundary of ${\bf{X}}$, which is not assumed to be smooth, and  assume that for some sufficiently small $\rho >0$, $\vol ( \gd {\bf{X}})_\rho \leq C \rho$, where $(\gd{\bf{X}})_\rho=\mklm{ x \in \R^n: \dist (x, \gd {\bf{X}}) <\rho}$.

Since $A$ commutes with the action of $G$ due to the invariance of $a$, the eigenspaces of $A$ are unitary $G$-modules  that decompose into irreducible subspaces.  In 1972,  Arnol'd \cite{arnold03} conjectured  that by studying the asymptotic behaviour of the spectral counting function
\bqn
N_\chi(\lambda)=d_\chi \sum_{t\leq \lambda} \mu_\chi(t) 
\eqn
where $\mu_\chi(\lambda)$ is the multiplicity of the irreducible representation of dimension $d_\chi$ corresponding to the character $\chi$ in the eigenspace of $A$ with eigenvalue $\lambda$,  one should be able to show  that the multiplicity of each unitary irreducible representation in  the above decomposition of $\L^2({\bf{X}})$ is asymptotically proportional to its dimension. 

The asymptotic distribution of eigenvalues was first studied by Weyl \cite{weyl} for certain second order differential operators in $\R^n$ using variational techniques. Another approach, which also gives an asymp\-totic description for the eigenfunctions,  was introduced by Carleman \cite{carleman}. His idea was to study the kernel of the resolvent, combined with a Tauberian argument. Minakshishundaram and  Pleijel \cite{minakshisundaram-pleijel} showed that one can study the Laplace transform of the spectral function as well, and extended the results of Weyl to closed manifolds, and G{\aa}rding \cite{garding} generalized Carleman's approach to higher order elliptic operators on bounded sets in $\R^n$. H\"ormander \cite{hoermander68} then extended these results to elliptic differential  operators on closed manifolds using the theory of Fourier integral operators. Further developments in this direction were given by Duistermaat and Guillemin, Helffer and Robert, and Ivrii. The first ones to study Weyl asymptotics for elliptic operators on closed Riemannian manifolds  in the presence of a compact group of isometries in a systematic way were Donnelly \cite{donnelly} together with  Br\"uning and Heintze \cite{bruening-heintze}, giving first order Weyl asymptotics for the spectral distribution function for each of the isotypic components, together with an estimate for the remainder in some special cases. Later, Guillemin and Uribe \cite{guillemin-uribe} described the relation between the spectrum of the considered operators, and the reduction of the corresponding bicharacteristic flow, and Helffer and Robert \cite{helffer-robert84,helffer-robert86} studied the situation in $\R^n$. Our approach is based on the Weyl calculus of pseudodifferential operators  developed by H\"ormander \cite{hoermander79}, and the method of approximate spectral projections, first introduced by Tulovskii and Shubin \cite{tulovsky-shubin}.  This  method is  somehow more closely related to the  original work of Weyl,  and starts from the observation that the asymptotic distribution function $N(\lambda)$ for the eigenvalues of an elliptic, self-adjoint operator is given by the trace of the orthogonal projection  on the space spanned by the eigenvectors corresponding to eigenvalues $\leq \lambda$. By introducing suitable approximations to these spectral projections in terms of  pseudodifferential operators, one can then derive asymptotics for $N(\lambda)$, and also obtain  estimates for the remainder term. Nevertheless, due to the presence of the boundary, the original method of Shubin and Tulovskii cannot be applied to our situation, and one is forced to use more elaborate techniques, which were subsequently developed by Feigin \cite{feigin} and Levendorskii \cite{levendorskii}. Recently, Bronstein and Ivrii have obtained even sharp estimates for the remainder term in the case of differential operators on manifolds with boundaries satisfying the conditions specified above.

This paper is structured as follows. Part I provides the foundations of the calculus of approximate spectral projection operators, and addresses the case where $G$ is a finite group of isometries. The case of a compact group of isometries will be the subject of Part II. The main result  of Part I is the following
 \begin{theorem}
 \label{thm:1}
Let $G$ be a finite group of isometries. Then the spectrum of $A$ is discrete, and the number  $N_\chi(\lambda)$ of  eigenvalues of $A$, counting multiplicities, less or equal $\lambda$ and with eigenfunctions in the  $\chi$-isotypic component $\mathrm{res} \, \H_\chi$ of $\L^2({\bf{X}})$, is given by 
\begin{equation*}
N_\chi(\lambda)=d_\chi  \sum_{t\leq \lambda} \mu_\chi(t) =\frac {d^2_\chi}{|G|}\gamma   \lambda^{n/2m} +O(\lambda^{(n-\epsilon)/2m})
\end{equation*}
for arbitrary  $\epsilon \in (0,\frac 1 2)$, where $|G|$ denotes the cardinality of $G$, $d_\chi$  the dimension of the irreducible representation of $G$ corresponding to the character $\chi$, and 
\bdm
\gamma=\frac 1 {n(2\pi)^n } \int_{\bf{X}} \int_{S^{n-1}} (a_{2m}(x,\xi))^{-n/2m} dx \, d \xi.  
\edm
Consequently, the multiplicity in $\L^2({\bf{X}})$ of the  irreducible representation corresponding to the character $\chi$  is given asymptotically by $\frac {d_\chi}{|G|} \gamma \lambda^{n/2m}$ as $\lambda \to \infty$.
\end{theorem}

\medskip

ACKNOWLEDGMENTS. The author wishes to thank Professor Mikhail Shubin for introducing him to this subject, and for many helpful discussions and useful remarks.

 \section{The Weyl Calculus for Pseudodifferential Operators in $\R^n$}

We first introduce the relevant symbol classes, as defined  in \cite{hoermander79}, and recall some theorems of Weyl calculus that will be needed in the sequel. We then study the pullback of symbols, and the composition of pseudodifferential operators with linear transformations. Thus, let $g$ be a slowly varying Riemannian metric in $\R^{l}$, regarded as a positive definite quadratic form,  and assume  that $m$ is a positive, $g$-continuous function on $\R^{l}$ (see  Definitions 2.1 and 2.2 in  \cite{hoermander79}).
 \begin{definition}
The class of symbols $S(g,m)$ is defined as the set of all functions $u\in \Cinft(\R^{l})$ such that, for every integer $k \geq 0$,
\bqn
\nu_k(g,m; u)= \sup_{x\in \R^{l}} \sup_{t_j \in \R^{l}} |u^{(k)}(x; t_1,\dots,t_k)|\Big /\Big ( \prod_{j=1}^k g_x(t_j)^{1/2} m(x)\Big )   < \infty.
\eqn
\end{definition}
Here $u^{(k)}$ stands for the $k$-th differential of $u$. 
Note that with the topology defined by the above semi-norms, $S(g,m)$ becomes a Fr\'{e}chet space. Consider now $\R^l = \R^n\oplus \R^n$,  regarded as  a symplectic space with the symplectic form
\bqn
\sigma(x,\xi;y,\eta)=\langle \xi,y \rangle -\langle x,\eta \rangle,
\eqn
where $\langle \cdot, \cdot \rangle$ denotes the usual Euclidean product of two vectors. Thus, $\sigma=\sum d\xi_j \wedge dx_j$. Assume that $g$ is $\sigma$-temperate, and that $m$ is $\sigma,g$-temperate (see Definition 4.1 in \cite{hoermander79}).  
If $a \in S(g,m)$ is interpreted as a Weyl symbol, the corresponding  pseudodifferential operator is given by
\bqn
\Op^w(a) u(x)= \int \int e^{i(x-y)\xi} a\Big (\frac {x+y} 2, \xi \Big )u(y) dy \, \dbar \xi,
\eqn
where $\dbar \xi= (2\pi)^{-n} d\xi$. Here and it what follows, it will be understood that each integral is to be performed over $\R^n$, unless otherwise specified. According to \cite{hoermander79}, Theorem 5.2, $\Op^w(a)$ defines a continuous linear map from $\S(\R^n)$ to $\S(\R^n)$, and from $\S'(\R^n)$ to $\S'(\R^n)$, and the corresponding class of operators will be denoted by $\Lcal(g,m)$. 
Moreover, one has the the following result concerning the   $\L^2$-continuity of pseudodifferential operators. 
\begin{theorem}
\label{thm:L2}
Let $g$ be a $\sigma$-temperate metric in $\R^n \oplus \R^n$, $g^\sigma$ the dual metric to $g$ with respect to $\sigma$, and $g \leq g^\sigma$.
Let  $a \in S(g,m)$, and assume that  $m$ is  $\sigma$, $g$-temperate.  Then $\Op^w(a):\L^2(\R^n) \rightarrow \L^2(\R^n)$ is a continuous operator if, and only if, $m$ is bounded. \end{theorem}
\begin{proof}
See \cite{hoermander79}, Theorem 5.3.
\end{proof}
The composition of pseudodifferential operators is described by the Main Theorem of Weyl Calculus.
\begin{theorem}
\label{thm:MT}
Let $g$ be a $\sigma$-temperate metric in $\R^n \oplus \R^n$,  and $g \leq g^\sigma$. Assume that  $a_1 \in S(g,m_1)$, $ a_2 \in S(g,m_2)$, where $m_1,\, m_2$ are $\sigma,\, g$-temperate functions. Then the composition of $\Op^w(a_1)$ with $\Op^w(a_2)$ in each of the spaces $\S(\R^n)$ or $\S'(\R^n)$ is a pseudodifferential operator with Weyl symbol $\sigma^w(\Op^w(a_1) \Op^w(a_2))$ in the class $S(g,m_1m_2)$. Moreover,
\begin{align}
\label{eq:A}
\begin{split}
\sigma^w(\Op^w(a_1) \Op^w(a_2))(x,\xi)&-\sum _{j<N} \big (\frac 12 i \sigma(D_x,D_\xi;D_y,D_\eta)\big )^j a_1(x,\xi) a_2(y,\eta)_{|(y,\eta)=(x,\xi)}/j!\\ & \in S(g, h_\sigma^N m_1m_2)
\end{split}
\end{align}
for every integer $N$, where $D_j=-i\gd_j$, $D=(D_1,\dots,D_n)$, and 
\bqn
h^2_\sigma(x,\xi)=\sup_{y,\eta} \frac {g_{x,\xi}(y,\eta)}{g_{x,\xi}^\sigma(y,\eta)}.
\eqn
\end{theorem}
\begin{proof}
See \cite{hoermander79}, Theorems 4.2 and 5.2.
\end{proof}
Note that  $g$ can always be written in the form $g(y,\eta)=\sum (\lambda_j y_j^2+ \mu_j \eta_j^2)$. Then $g^\sigma(x,\xi)=\sum  (y_j^2/\lambda_j+ \eta_j^2/\mu_j)$, so that 
\bq
\label{eq:2a}
h_\sigma(x,\xi)=\max(\lambda_j \mu_j)^{1/2}.
\eq
The following proposition describes the asymptotic expansion of symbols, see \cite{levendorskii}, Theorem 3.3.
\begin{proposition}
\label{prop:A}
Let $a_j \in S(g, h_\sigma ^{N_j} m)$ be a sequence of symbols such that $0=N_1<N_2< \dots \rightarrow \infty$. Then there exists a symbol $a \in S(g,m)$ such that 
\begin{itemize}
\item[a)] $\supp a \subset \bigcup_{j} \supp a_j$;
\item[b)] $a - \sum_{j=1}^{l-1} a_j \in S(g,h_\sigma^{N_l} m), \quad l >1$.
\end{itemize}
In this case, one writes $a \sim \sum a_j$.
\end{proposition}
We will further write 
\bdm
S^{-\infty}(g,m)= \bigcap _{N=1}^\infty S(g,h_\sigma^Nm),
\edm
and denote the corresponding operator class by $\Lcal^{-\infty}(g,m)$.
We introduce now certain  hypoelliptic symbols which will be needed in the sequel. They were introduced by Levendorskii in \cite{levendorskii}.
\begin{definition}
The class of symbols $SI(g,m)$ consists of all $a \in S(g,m)$ that can be represented in the form $a=a_1+a_2$, where $c m < |a_1|$ and $a_2 \in S(g,h_\sigma^\epsilon m)$ for some constants $c,\epsilon >0$. The corresponding class of operators is denoted by $\mathcal{LI}(g,m)$. If instead of $cm < |a_1|$ one has $cm<a_1$, one  writes $a \in SI^+(g,m)$ and $\mathcal{LI}^+(g,m)$, respectively. 
\end{definition}
For a proof of the following lemmas, we refer the reader to \cite{levendorskii}, Lemma 5.5, Lemma 8.1, and Lemma 8.2.
\begin{lemma}
\label{lemma:SI1}
Let $a \in SI(g,m)$. Then there exists a symbol $b \in SI(g,m^{-1})$ such that 
\bqn
\Op^w(a) \Op^w(b) -\1 \in \Lcal^{-\infty} (g,1), \qquad Op^w(b) \Op^w(a) -\1 \in \Lcal^{-\infty} (g,1).
\eqn
The operator $\Op^w(b)$ is called a parametrix for $\Op^w(a)$.
\end{lemma}
\begin{lemma}
\label{lemma:SI2}
If $a \in SI^+(g,m)$, then there exists a symbol $b \in S(g,m^{1/2})$ such that
\bqn
\Op^w(a)-\Op^w(b)^\ast  \Op^w(b) \in \Lcal^{-\infty}(g, m),
\eqn 
where $ \Op^w(b)^\ast$ is the adjoint of $ \Op^w(b)$.
 \end{lemma}
\begin{lemma}
\label{lemma:3}
Let $\epsilon >0$, and $a_t\in S(g,h_\sigma^\epsilon)$, $t \in \R$,  be a family of symbols depending on a parameter. Furthermore, assume that the corresponding seminorms $\nu_k(g,h_\sigma^\epsilon;a_t)$ are bounded by some constants independent of $t$, and let $c >0$ be arbitrary. Then there exists a subspace $L \subset \L^2(\R^n)$ of finite codimension such that 
\bqn
\norm{\Op^w(a_t)u}_{\L^2} \leq c \norm{u}_{\L^2} \qquad \text{for all } u \in L  \text{ and all } t\in \R.
\eqn  
\end{lemma}
\begin{remark}
\label{rem:L2}
Lemma \ref{lemma:3} is a consequence of the fact that, for $a \in S(g,1)$, one has the uniform bound
\bqn
\norm{\Op^w(a)}_{\L^2} \leq C \max _{k\leq N} \, \nu_k(g,1;a),
\eqn
where  $C>0$ and $N\in \N$ depend only on the constants characterizing $g$, but not on $a$ (see the proof of the sufficiency in Theorem 5.3 in \cite{hoermander79}, and Theorem 4.2 in \cite{levendorskii}).
\end{remark}
In general, the pullback of symbols under $\Cinft$ mappings is described by the following
\begin{lemma}
\label{lemma:5a}
Let $g_1,g_2$ be slowly varying metrics on $\R^l$, respectively $\R^{l'}$, and $\chi \in \Cinft (\R^l,\R^{l'})$. Then 
\bqn
\chi^\ast S(g_2,1) \subset S(g_1,1)
\eqn
if, and only if, for every $k>0$,
\bqn
g_{2\chi(x)} ( \chi^{(k)}(x; t_1,\dots,t_k)) \leq C_k \prod _{j=1}^k g_{1x}(t_j), \qquad x, t_1,\dots,t_k \in \R^l.
\eqn
In particular, if $m$ is $g_2$-continuous, then $\chi^\ast m$ is $g_1$-continuous and $\chi^\ast S(g_2,m) \subset S(g_1,\chi^\ast m)$.
\end{lemma}  
\begin{proof}
See \cite{hoermander79}, Lemma 8.1.
\end{proof}
In all our applications, we will be  dealing mainly with metrics $g$ on $\R^{2n}$ of the form
\bq
\label{eq:6}
g_{x,\xi}(y,\eta)= (1+|x|^2 +|\xi|^2)^\delta |y|^2+ (1+|x|^2+ |\xi|^2)^{-\rho} |\eta|^2,
\eq
 where  $1 \geq \rho>\delta \geq0$. The conditions of Theorem \ref{thm:MT} are satisfied then, and $h^2_\sigma(x,\xi)=(1+|x|^2+|\xi|^2)^{{\delta -\rho}}$ by \eqref{eq:2a}. For the rest of this section, assume that  $g$ is  of the form \eqref{eq:6}, and put $h(x,\xi)=(1+|x|^2+|\xi|^2)^{-1/2}$. In this case, the space of symbols $S(g,m)$ can also be characterized as follows.
 
   \begin{definition} Let $g$ be of the form \eqref{eq:6}, and $m$ a $g$-continuous function.
  The class $\Gamma_{\rho,\delta}(m, \R^{2n})$, $0\leq \delta < \rho \leq 1$, consists of all functions $u \in \Cinft(\R^{2n})$ which for all multiindices $\alpha,\beta$ satisfy the estimates
 \bqn
 |\gd ^\alpha _\xi \gd^\beta_x u(x,\xi)| \leq C_{\alpha \beta}\,  m(x,\xi) \, (1+|x|^2+|\xi|^2)^{(-\rho|\alpha|+\delta|\beta|)/2}.
  \eqn
    In particular, we will  write  $\Gamma^l_{\rho,\delta}(\R^{2n})$ for $\Gamma_{\rho,\delta}(h^{-l},\R^{2n})$, where  $l \in \R$.
  \end{definition}
  An easy computation then shows that  $S(g,m)=\Gamma_{\rho,\delta}(m,\R^{2n})$.
  For future reference, note that $u \in S(g,m)$ implies $\gd ^\alpha_\xi \gd^\beta_x u \in S(g, m \, h^{\rho|\alpha|- \delta|\beta|})$. The pullback of symbols for metrics of the form \eqref{eq:6} can now be described as follows.
 \begin{lemma} 
 \label{lemma:5}
  Let  $ \delta +\rho \geq 1$, and $g$ be a metric of the form \eqref{eq:6}. Assume that   $\chi(x,\xi)=(y(x),\eta(x,\xi))$ is  a diffeomorphism in $\R^{2n}$ such that $\eta$ is linear in $\xi$, and the derivatives of $y$ and $\eta$ are bounded for $|\xi|<1$. Furthermore, let
 \bqn 
 \frac 1 C g_{x,\xi} (t) \leq g_{\chi(x,\xi)} (t) \leq C g_{x,\xi}(t), \qquad \frac 1 C m(x,\xi)\leq  \chi^\ast m(x,\xi) \leq C m(x,\xi),
 \eqn
 where $m$ is a $g$-continuous function, and $C>0$ is a suitable constant. Then $\chi^\ast S(g,m) \subset S(g,\chi^\ast m)$.
  \end{lemma}
  \begin{proof} Instead of verifying the necessary and sufficient condition in Lemma \ref{lemma:5a}, we will prove the statement directly.
Let $b \in S(g,m)$, and let $s,t\dots$ be $k$ vectors in $\R^{2n}$. The $k$-th differential 
\bqn
(b\circ \chi)^{(k)}(x,\xi;s,t\dots)=\langle t,D \rangle \langle s,D\rangle \dots (b\circ \chi) (x,\xi)
\eqn
 is given by a sum of terms of the form $s_i t_j \dots \gd^\alpha (b\circ \chi )(x,\xi)$, where we can assume that all the coefficients $s_i, t_j,\dots$ are different from zero; in particular, $(b\circ \chi)^{(1)}(x,\xi)=b^{(1)}(\chi(x,\xi)) \chi^{(1)}(x,\xi)$, where 
 \bqn
 \chi^{(1)} (x,\xi)=\left ( \begin{array}{cc} y^{(1)}(x) & 0 \\ A(x,\xi) & B(x) \end{array} \right ),
 \eqn
 $A$ being linear in $\xi$. The derivatives $\gd ^\alpha (b\circ \chi)(x,\xi)$ are sums of expressions of the form
 $$(\gd ^\beta b)(\chi(x,\xi))(\gd^{\gamma_1}\chi_{i_1}) (x,\xi)\dots (\gd^{\gamma_{l}}\chi_{i_l}) (x,\xi),$$
 where $\gamma_1 + \dots + \gamma_l=\alpha$ and $l=|\beta|$. Since additional powers of $\xi$ only appear in companion with additional derivatives of $b$ with respect to $\eta$ that originate from derivatives of $b \circ \chi$ with respect to $x$, each of the terms of $(b\circ \chi)^{(k)}(x,\xi;s,t, \dots)$ can  be estimated from above by some constant times an expression of the form
 \bqn
  \,  | s_it_j\dots (\gd_y^{\beta'}\gd_\eta^{\beta''}b)(\chi(x,\xi)) P^d(x,\xi)|,
 \eqn
 where $P^d(x, \xi)$ is a homogeneous polynomial in $\xi$ of degree $d$ which is bounded for $|\xi|<1$, and  
 \bqn
d= |\beta''|-N''=|\beta''|-k+N';
 \eqn
  here $N'=|\alpha'|$ and $N''=|\alpha''|$ denote the number of $x$- and $\xi$-components in the product $s_it_j\dots$, respectively. Indeed, if we differentiate in $\gd ^\alpha (b \circ \chi) (x,\xi)$ first with respect to $\xi$ we get
\bqn
\gd ^{\alpha''}_\xi (b \circ \chi)(x,\xi)=\sum_{\eta_{j_1}, \dots, \eta_{j_{|\alpha''|}}=1}^n (\gd_{\eta_{j_1}} \dots \gd_{\eta_{j_{|\alpha''|}}} b)(\chi(x,\xi))  \frac {\gd \eta_{j_1}}{\gd \xi_{i_1}}(x) \dots \frac {\gd \eta_{j_{|\alpha''|}}}{\gd \xi_{i_{|\alpha''|}}}(x),
\eqn
where $\gd_\xi^{\alpha''}=\gd_{\xi_{i_1}} \dots \gd_{\xi_{i_{|\alpha''|}}}$,
and differentiating now with respect to $x$ yields the assertion. Note that $N''\leq |\beta''|$.
In order to prove the assertion of the lemma,  we have to show that
 \bq
 \label{3a}
 \sup_{x,\xi} \sup_{s,t,\dots} \frac {| s_it_j\dots (\gd_y^{\beta'}\gd_\eta^{\beta''}b)(\chi(x,\xi)) P^d(x,\xi)|}{g_{x,\xi}^{1/2}(s) g_{x,\xi}^{1/2} (t) \dots m(x,\xi)}<\infty,
 \eq
where it suffices to consider only those $s,t, \dots$ whose only non-zero components are $s_i,t_j \dots$.
 Since $N' \geq d$,  there are $d$ vectors $p,q,\dots$ among the vectors $s,t,\dots$ contributing with $x$-components to the product $s_it_jp_kq_l \dots$. Furthermore, let $w,z,\dots$ be $d$ vectors such that $w_{n+k}=p_k, z_{n+l}=q_l, \dots$, their other components being zero. We then obtain the estimate
 \begin{gather*}
 \frac {|s_it_jp_kq_l \dots (\gd_x^{\beta'}\gd_\xi^{\beta''}b)(\chi(x,\xi))|}{m(\chi(x,\xi))g_{\chi(x,\xi)}^{1/2}(s) g_{\chi(x,\xi)}^{1/2}(t)\dots g_{\chi(x,\xi)}^{1/2}(w)g_{\chi(x,\xi)}^{1/2}(z)\dots } \cdot \frac {| P^d(\xi) g_{x,\xi}^{1/2}(w)g_{x,\xi}^{1/2}(z)\dots |}{g_{x,\xi}^{1/2}(p)g_{x,\xi}^{1/2}(q)\dots } \\
 \leq C (1 + |x|^2 + |\xi|^2)^{d(1-\delta-\rho)/ 2 } 
 \end{gather*}
for all $x,\xi, s, t, \dots$. Indeed,  
 \bqn
 g_{x,\xi}^{1/2}(w) = |p_k| (1 + |x|^2 +|\xi|^2 ) ^{-\rho/2}, \qquad g_{x,\xi}^{1/2}(p) =|p_k| (1 + |x|^2 +|\xi|^2 ) ^{\delta/2},\dots.
 \eqn
 On the other hand, besides the $d$ vectors $p,q\dots$ there are still $N'-d=k-|\beta''| \geq |\beta'|$ vectors among the remaining vectors $s,t\dots$ contributing with $x$-components to the product $s_it_j\dots$. Since the corresponding quotients
$
{|r_l|}/{g^{1/2}_{\chi(x,\xi)}(r)}$ 
can be estimated from above by some constant, we can assume that there are precisely $|\beta'|$ of them. Also note that there are exactly $d +N''=|\beta''|$ vectors among the vectors $s,t\dots w,z\dots$ contributing with $\xi$-components to $s_it_j\dots$. We can therefore assume that the components of $s,t\dots w,z,\dots$ are prescribed by the multiindex $\beta=(\beta',\beta'')$ in such a way that
$$s_i t_jp_kq_l \dots (\gd_x^{\beta'} \gd _\xi^{\beta''} b) ( \chi(x,\xi)) = b^{(|\beta|)} (\chi(x,\xi);s,t\dots w,z\dots).$$
The desired estimate \eqref{3a}  now follows by using the assumptions that $b \in S(g,m)$ and $\delta+ \rho\geq 1$.
 \end{proof}
  If $a\in S(g,m)$ is regarded as a right, respectively left symbol, the corresponding pseudodifferential operators are given by
\bqn
\Op^l(a) u(x)= \int \int e^{i(x-y)\xi} a (x, \xi  )u(y) dy \, \dbar \xi, \quad \Op^r(a) u(x)= \int \int e^{i(x-y)\xi} a (y, \xi  )u(y) dy \, \dbar \xi,
\eqn
where $g$ is assumed to be of the form \eqref{eq:6}. 
By \cite{hoermander79}, Theorem 4.5, the three sets of operators $\Op^w(a)$, $\Op^l(a)$, and $\Op^r(a)$ coincide.
Theorem \ref{thm:MT} can then also be formulated in terms of left and right symbols.
  In what follows, we would like to treat left, right, and Weyl symbols on the same grounding by introducing the notion of the $\tau$-symbol. To do so, we introduce yet another class of amplitudes which is closely related to the space $\Gamma^l_{\rho,\delta}(\R^{2n})$, compare \cite{shubin}, Chapter 4. 
  \begin{definition}
  The class $\Pi ^l_{\rho,\delta}(\R^{3n})$ consists of all functions $u \in \Cinft(\R^{3n})$ which for a suitable $l'\in \R$ satisfy the estimates
  \bqn
|\gd_\xi^\alpha \gd_x^\beta \gd_y^\gamma u (x,y,\xi)| \leq C_{\alpha\beta\gamma} (1+|x|^2+|y|^2+|\xi|^2)^{(l-\rho|\alpha|+\delta|\beta+\gamma|)/2} (1+|x-y|^2)^{(l'+\rho|\alpha|+\delta|\beta+\gamma|)/2}.
  \eqn
  \end{definition}
  The relationship between the spaces $\Pi ^l_{\rho,\delta}(\R^{3n})$ and $\Gamma^l_{\rho,\delta}(\R^{2n})$ is described by the following lemma.
  \begin{lemma}
  \label{lemma:6}
  Let $0 \leq \delta < \rho \leq 1$, and $p: \R^{2n} \rightarrow \R^n$ be a linear map such that $(x,y) \mapsto (p(x,y),x-y)$ is an isomorphism. Let $a (w,\eta) \in \Gamma^l_{\rho,\delta}(\R^{2n})$, and define
  \bqn
  b(x,y,\xi) = a(p(x,y), \psi (x,y) \xi),
  \eqn
  where $\psi: \Xi \rightarrow \GL(n,\R)$ is a $\Cinft$ mapping on some open subset $\Xi \subset \R^{2n}$, having  bounded derivatives. If $\delta+\rho \geq 1$, then $b \in \Pi ^l_{\rho,\delta}(\Xi \times \R^{n})$.
 \end{lemma}
  \begin{proof}
  We will proof the assertion by induction on $|\alpha+\beta+\gamma|$. First note that $\gd^\alpha_\xi \gd ^\beta _x \gd ^\gamma_y b(x,y,\xi)$ is given by a sum of terms of the form
  \bq
  \label{eq:7}
  (\gd^{\alpha'}_\eta\gd^{\beta'}_w a)(p(x,y),\psi(x,y) \xi) P^d(x,y,\xi),
  \eq
  where $P^d(x,y,\xi)$ is a polynomial in $\xi$ of degree $d$. Each of these summands can  be estimated from above by
  \bqn
  C ( 1 +|p(x,y)|^2 +|\psi(x,y) \xi|^2)^{(l-\rho|\alpha'| +\delta|\beta'|)/2} |P^d(\xi)|,
  \eqn
where $P^d(\xi)$ is a polynomial in $\xi$ of degree $d$ with constant coefficients, and $C>0$ is a constant. We assert that  the inequality 
  \bq
  \label{eq:8}
  -\rho|\alpha'| +\delta |\beta'| +d \leq -\rho |\alpha| + \delta |\beta +\gamma|
  \eq
  holds for all $|\alpha+\beta+\gamma|=N$, and all occurring combinations of $\alpha'$, $\beta'$, and $d$. It is not difficult to verify the assertion for $N=1$. Let us now assume that  \eqref{eq:8} holds for $|\alpha+\beta+\gamma|=N$. Differentiating \eqref{eq:7} with respect to $\xi_j$ yields
  \begin{align*}
  \sum _{i=1}^n &(\gd_{\eta_i} \gd_{\eta}^{\alpha'}\gd_w^{\beta'}a)(p(x,y),\psi(x,y) \xi) \psi(x,y)_{ij} P^d(x,y,\xi)\\& +(\gd^{\alpha'}_\eta\gd^{\beta'}_w a)(p(x,y),\psi(x,y) \xi) \gd_{\xi_j}P^d(x,y,\xi),
  \end{align*}
  and we get the inequalities
  \begin{align*}
  \begin{split}
 & -\rho(|\alpha'|+1)+\delta|\beta'| +d \leq -\rho (|\alpha|+1) +\delta |\beta+\gamma|, \\
  &-\rho |\alpha'|+\delta|\beta'|+d-1  \leq -\rho(|\alpha|+1) +\delta |\beta+\gamma|.
  \end{split}
  \end{align*}
  Similarly, differentiation with respect to, say $x_j$, gives
   \begin{align*}
   \sum _{i=1}^n &(\gd_{w_i} \gd_{\eta}^{\alpha'}\gd_w^{\beta'}a)(p(x,y),\psi(x,y) \xi) (\gd_{x_j} p_i)(x,y)  P^d(x,y,\xi)\\
 &+\sum _{i=1}^n (\gd_{\eta_i} \gd_{\eta}^{\alpha'}\gd_w^{\beta'}a)(p(x,y),\psi(x,y) \xi) \gd_{x_j} (\psi(x,y) \xi)_{i} P^d(x,y,\xi)\\ 
 &+(\gd^{\alpha'}_\eta\gd^{\beta'}_w a)(p(x,y),\psi(x,y) \xi) (\gd_{x_j}P^d)(x,y,\xi),
  \end{align*}
  and we arrive at  the inequalities
  \begin{align*}
  \begin{split}
  &-\rho|\alpha'| +\delta(|\beta'|+1) +d \leq -\rho|\alpha| + \delta(|\beta+\gamma|+1),\\
   & -\rho(|\alpha'|+1)+\delta|\beta'| +d+1 \leq -\rho |\alpha|+\delta |\beta+\gamma|-\rho+1 \leq -\rho|\alpha|+\delta (|\beta+\gamma|+1), \\
  &-\rho |\alpha'|+\delta|\beta'|+d \leq -\rho|\alpha|+\delta (|\beta+\gamma|+1),
  \end{split}
  \end{align*}
where, in particular, we made use of the assumption $\delta+\rho \geq 1$. This proves \eqref{eq:8} for $|\alpha+\beta+\gamma|=N+1$. Summing up, we get the estimate
\begin{align*}
|\gd^\alpha_\xi \gd ^\beta _x \gd ^\gamma_y b(x,y,\xi)|&\leq C_1(1+|p(x,y)|+|\xi|)^{l-\rho|\alpha|+\delta|\beta+\gamma|}\\
&\leq C_2 (1+|(p(x,y)| +|x-y| +|\xi|)^{l-\rho |\alpha| +\delta |\beta+\gamma|}(1+|x-y|)^{|l|+\rho|\alpha|+\delta|\beta+\gamma|},
\end{align*}
 where the latter inequality follows by using the easily verified inequality
\bqn
\frac {(1+ |p(x,y)|+|\xi|)^s}{(1+|p(x,y)|+|x-y|+|\xi|)^s} \leq C (1+|x-y|)^{|s|}, \qquad s \in \R,
\eqn
compare the proof of Proposition 23.3 in \cite{shubin}. Since $|x|+|y|$ and $|p(x,y)|+|x-y|$ define equivalent metrics, the assertion of the lemma follows.  
  \end{proof}
  \begin{proposition}
  Let $a (x,y,\xi) \in \Pi^l_{\rho,\delta}(\R^{3n})$, where $1\geq \rho > \delta \geq 0$. Then the oscillatory integral
  \bq
  \label{eq:OI}
  Au(x)=\int \int e^{i(x-y)\xi} a(x,y,\xi) u(y) dy \, \dbar \xi
  \eq
  defines a continuous linear operator from $S(\R^n)$ to $\S(\R^n)$, and from $\S'(\R^n)$ to $\S'(\R^n)$.
  \end{proposition}
  \begin{proof}
  Consider first the case $a \in \CT(\R^{3n})$, and assume that $u \in \Cinft(\R^n)$ has bounded derivatives. Then the integration in \eqref{eq:OI} is carried out over a compact set, and partial integration gives
  \bqn
  Au(x)= \int \int e^{i(x-y)\xi} \eklm{x-y}^{-M} \eklm{D_\xi}^M \eklm{D_y}^N [\eklm{\xi}^{-N} a (x,y,\xi) u(y)] dy \, \dbar \xi,
  \eqn 
  where $M, N$ are even non-negative integers, and  $\eklm{x}$ stands for $(1+x_1^2+\dots+x_n^2)^{1/2}$. 
  Let now $a \in \Pi^l _{\rho,\delta} (\R^{3n})$, and assume that $M,N$ are such that $l-N(1-\delta) <-n$, $ l + l'+2\delta N-M< -n$.   The latter integral then becomes absolutely convergent, defining a continuous function of $x$, and represents the regularization of the oscillatory integral \eqref{eq:OI}. Increasing $M$ and $N$ we will obtain integrals which are convergent also after differentiation with respect to $x$. In view of the inequality $\eklm{x} ^k \leq \eklm{y}^k \eklm{x-y}^k$, where $k>0$,  one finally sees that  $A$ defines a continuous map from $\S(\R^n)$ to $\S(\R^n)$, which, by duality, can be extended to a continuous map from $\S'(\R^n)$ to $\S'(\R^n)$.  
  \end{proof}
  We can now introduce the notion of the $\tau$-symbol. In what follows, $m$ will be a $g$-continuous function.
  \begin{corollary}
Let $a \in S(g,m)=\Gamma_{\rho,\delta}(m,\R^{2n})$,  $0 \leq 1-\rho \leq \delta< \rho \leq 1$, and $\tau \in \R$. Then
  \bqn
Au(x) = \int \int e^{i(x-y)\xi} a ((1-\tau)x+\tau y, \xi  )u(y) dy \, \dbar \xi
  \eqn
  defines a continuous operator in $\S(\R^n)$, respectively $\S'(\R^n)$. In this case, $a$ is called the $\tau$-symbol of $A$, and the operator $A$ is denoted by  $\Op^\tau(a)$.
  \end{corollary}
 \begin{proof}
For simplicity, we restrict ourselves to the  case  $m=h^{-l}$. By Lemma \ref{lemma:6} we then have $b(x,y,\xi)=a((1-\tau)x+\tau y,\xi) \in \Pi^l_{\rho,\delta}(\R^{3n})$, and the assertion follows with the previous proposition. The case of a general  $m$ is proved in a similar way.
 \end{proof}
 Our next aim is to prove the following
 \begin{theorem}
 \label{thm:3}
 Let $0 \leq 1-\rho \leq \delta< \rho \leq 1$,  $\tau,\tau' \in \R$ be arbitrary, $a(x,\xi) \in S(g,m)=\Gamma_{\rho,\delta}(m,\R^{2n})$,  and assume that $\kappa:\R^n\rightarrow \R^n$ is an invertible linear map. Furthermore, assume that $m$ is invariant under $\kappa$ in the sense that $ m(\kappa^{-1}(x), \,^t \kappa\,  (\xi)) = m(x,\xi)$, 
  and set $A=\Op^{\tau'} (a)$. Then 
 \bqn
 A_1u=[A(u\circ \kappa)] \circ \kappa^{-1}, \qquad u \in \S(\R^n),
 \eqn
 defines a pseudodifferential operator with a uniquely defined $\tau$-symbol $\sigma^{\tau}(A_1) \in S(g,m)$.
 \end{theorem} 
  \begin{proof}
 Let us consider first  the  case $m=h^{-l}$.
  Putting  $\kappa_1=\kappa^{-1}$, one sees that  $A_1$ is a Fourier integral operator given by
  \begin{align*}
A_1u(x)&=\int \int e^{i(\kappa_1(x)-y)\cdot \xi} a((1-\tau')\kappa_1(x)+\tau' y,\xi) u (\kappa(y)) dy\,  \dbar \xi\\
&=\int \int e^{i(\kappa_1(x)-\kappa_1(y))\cdot  \xi} a((1-\tau')\kappa_1(x)+\tau' \kappa_1(y),\xi) |\det \kappa_1'(y)| u(y) dy \, \dbar \xi,
  \end{align*}
 and performing the change of variables $\xi \mapsto \, ^t \kappa\,(\xi)$, we get
\bqn
\label{eq:11}
A_1u(x) = \int \int e^{i(x-y)\cdot  \xi} a_1(x,y, \xi) u(y) dy \, \dbar \xi, 
\eqn
where we put $ a_1(x,y,\xi)= a((1-\tau')\kappa_1(x)+\tau'\kappa_1(y), ^t \kappa\,( \xi)) |\det \kappa_1| |\det \, ^t \kappa  |  $.  Applying  Lemma \ref{lemma:6} with $p(x,y)=(1-\tau') \kappa_1(x)+\tau'\kappa_1(y)$, one obtains $a_1(x,y,\xi) \in \Pi^l_{\rho,\delta}(\R^{3n})$ for arbitrary $a \in \Gamma^l_{\rho,\delta}(\R^{2n})=S(g,h^{-l})$. Next, let us introduce the coordinates $v=(1-\tau)x+\tau y$, $w=x-y$, and expand $a_1(x,y,\xi)=a_1(v+\tau w, v-(1-\tau)w,\xi)$  into a Taylor series at $w=0$, compare \cite{shubin}, pages 180-182. This yields
\bqn 
a_1(x,y,\xi) =\sum _{|\beta+\gamma| \leq N-1} \frac{(-1)^{|\gamma|}}{\beta! \gamma!}\tau^{|\beta|} (1 -\tau)^{|\gamma|} (x-y)^{\beta+\gamma} (\gd _x ^\beta \gd_y ^\gamma a_1)(v,v,\xi)+r_N(x,y,\xi),
\eqn
where
\bqn
r_N(x,y,\xi) = \sum_{|\beta+\gamma|=N} c_{\beta \gamma}(x-y)^{\beta+\gamma}\int^1_0 (1-t)^{N-1} (\gd_x^\beta \gd_y^\gamma a_1) (v+t\tau w, v-t(1-\tau)w,\xi) dt,
\eqn
$c_{\beta\gamma}$ being constants. 
Since  the operator with amplitude $(x-y)^{\beta+\gamma}(\gd_x^\beta\gd_y^\gamma a_1)(v,v,\xi) $ coincides with the one with amplitude $(-1)^{|\beta+\gamma|} ( \gd_\xi^{\beta +\gamma} D_x^\beta D_y^\gamma a_1) (v,v,\xi)$, we can write  $A_1$ also as $A_1=B_N+R_N$, where $B_N$ is the  operator with $\tau$-symbol
\bqn 
b_N(x,\xi)=\sum _{|\beta+\gamma| \leq N-1} \frac{1}{\beta! \gamma!}\tau^{|\beta|} (1 -\tau)^{|\gamma|}  \gd_\xi^{\beta +\gamma} (-D_x) ^\beta D_y ^\gamma a_1(x,y,\xi)_{|y=x},
\eqn
and $R_N$ has amplitude $r_N(x,y,\xi)$.  Similarly, we can assume that $R_N$ is given by a sum of terms having  amplitudes of  the form 
\bqn
\int_0 ^1 (\gd^{\beta+\gamma}_\xi \gd ^\beta _x \gd ^\gamma_y a_1)( v+t\tau w, v-t(1-\tau)w,\xi)(1-t)^{N-1} dt,
\eqn
where $|\beta +\gamma| =N$.
In view of the estimate
\bqn
|(\gd^{\beta+\gamma}_\xi \gd_x^\beta \gd^\gamma_y a_1) (v+t \tau w, v-t(1-\tau) w,\xi)| \leq C(1+|v| +|wt|+|\xi|)^{l-N(\rho-\delta)} (1+|tw|)^{l'+N(\rho+\delta) },
\eqn 
for some $l'$ and    $|\beta +\gamma|=N$, one can then show that  $r_N(x,y,\xi) \in \Pi_{\rho,\delta}^{l-N(\rho-\delta)}(\R^{3n})$, where, by assumption,  $\rho-\delta >0$. Define now $A_1'$ as the pseudodifferential operator with $\tau$-symbol 
\bq
\label{eq:11a}
a_1' (x,\xi) \sim \sum_{N=0}^\infty (b_N(x,\xi) -b_{N-1}(x,\xi)).
\eq
Then $A_1-A_1'$ has kernel and $\tau$-symbol belonging to $\S(\R^{2n})$. Since $b_N(x,\xi) \in S(g,h^{-l})$ for all $N$, the assertion of the theorem follows  in view of the uniqueness of the $\tau$-symbol, and $\sigma^{\tau}(A_1) \in \Gamma^l_{\rho,\delta}(\R^{2n})$. Let us consider now the case of a general $m$. By examining the proof of Lemma \ref{lemma:6}, we see that $a_1(x,y,\xi)$ must satisfy an estimate of the form
\begin{align*}
|\gd^\alpha_\xi \gd ^\beta _x \gd ^\gamma_y a_1(x,y,\xi)|&\leq C_1 \, m(p(x,y),\, ^t \kappa\, ( \xi))\, (1+|p(x,y)|+|\xi|)^{-\rho|\alpha|+\delta|\beta+\gamma|}\\
&\leq C_2 \,m(p(x,y),\, ^t \kappa\, (\xi))\,(1+|x|+|y| +|\xi|)^{-\rho |\alpha| +\delta |\beta+\gamma|}(1+|x-y|)^{\rho|\alpha|+\delta|\beta+\gamma|}.
\end{align*}
 Consequently,
\begin{gather*}
|(\gd^{\beta+\gamma}_\xi \gd_x^\beta \gd^\gamma_y a_1) (v+t \tau w, v-t(1-\tau) w,\xi)| \leq \\ C \,m(p(v+t\tau w, v-t(1-\tau)w),\, ^t \kappa\, ( \xi))\,(1+|v| +|wt|+|\xi|)^{-N(\rho-\delta)} (1+|tw|)^{l'+  N(\rho +\delta)},
\end{gather*}
where   $|\beta +\gamma|=N$, and we can  again define $A_1'=\mathrm{Op}^\tau(a_1')$  by the asymptotic expansion \eqref{eq:11a}, such that $A_1-A_1'$ has kernel and $\tau$-symbol belonging to $\S(\R^{2n})$. The  assertion of the theorem now  follows by noting that $b_N(x,\xi) \in S(g,m)=\Gamma_{\rho,\delta}(m,\R^{2n})$ for all $N$, due to the invariance of m.  In particular, one has the asymptotic expansion
\bq
\label{eq:12}
\sigma^{\tau}(A_1)(x,\xi) - \sum_{|\beta+\gamma|\leq N-1} \frac 1 {\beta! \gamma!} \tau^{|\beta|}(1-\tau)^{|\gamma|} \gd_\xi^{\beta+\gamma}(-D_x)^\beta D_y^\gamma a_1(x,y,\xi)_{|y=x} \in S(g, h_\sigma^N m)
\eq
for arbitrary integers $N$, where the first summand is given by $a_1(x,x,\xi)=a(\kappa^{-1}(x), \, ^t \kappa \, (\xi))$. 
\end{proof}
Theorem \ref{thm:3} allows us, in particular, to express the $\tau$-symbol of an operator in terms of its $\tau'$-symbol. More generally, one has the following 
\begin{corollary}
\label{cor:1}
In the setting of Theorem \ref{thm:3} assume that, in addition, $a(\kappa^{-1}(x),  \,^t \kappa \, (\xi))=a(x,\xi)$, and $\det \kappa =\pm 1$.  Then $A_1=A$, and the $\tau$-symbol  of $A=\Op^{\tau'}(a)$ is given by
\bq
\label{eq:12bis}
\sigma^{\tau}(A)(x,\xi) \sim \sum_{\beta,\gamma} \frac 1 {\beta! \gamma!} \tau^{|\beta|}(1-\tau)^{|\gamma|}(\tau'-1)^{|\beta|} {\tau'}^{|\gamma|} \gd_\xi^{\beta+\gamma}D_x^{\beta+\gamma}  a(x,\xi).
\eq
\end{corollary}
\begin{proof}
With $a_1(x,y,\xi)$ defined as in the proof of Theorem \ref{thm:3}, we have $a_1(x,y,\xi)=a((1-\tau')x +\tau'y,\xi)$, so that $A_1=\Op^{\tau'}(a)=A$. The corollary then follows with  the asymptotic expansion \eqref{eq:12}.
\end{proof}

 \section{The approximate spectral projection operators}

Let  $G\subset\O(n)$ be a compact group of isometries acting on Euclidean space $\R^n$, and ${\bf{X}}$ a bounded domain in $\R^n$ which is invariant under $G$.
Consider the regular representation $T$ in the Hilbert spaces  $\L^2(\R^n)$ and $\L^2({\bf{X}})$, respectively, and endow them with a $G$-invariant scalar product, so that  $T$ becomes unitary. Let $A_0$ be a symmetric, classical pseudodifferential operator of order $2m$ with principal symbol $a_{2m}$ as defined in \cite{shubin}, and regard it  as an operator in $\L^2(\R^n)$ with domain $\CT(\R^n)$. Furthermore, assume that $A_0$ is $G$-invariant, i.e. that it commutes with the operators $T(g)$ for all  $g\in G$, and that 
 \bq
 \label{eq:8a}
 (A_0u,u) \geq c \norm{u}_m^2, \qquad u \in \CT({\bf{X}}),
 \eq
for some $c>0$, where $(\cdot, \cdot)$ denotes the scalar product in $\L^2(\R^n)$, and $\norm{\cdot}_s$ is a norm in the Sobolev space $H^s(\R^n)$.  
Consider next   the decomposition of $\L^2(\R^n)$ and $\L^2({\bf{X}})$ into  isotypic components,
\bqn
\L^2(\R^n) =\bigoplus_{\chi \in \hat G} \H_\chi, \qquad \L^2({\bf{X}})= \bigoplus_{\chi \in \hat G} \mathrm{res} \, \H_\chi,
\eqn
where $\hat G$ is the set of all irreducible characters of $G$, and $\mathrm{res}$ denotes the restriction of functions defined on  $\R^n$ to ${\bf{X}}$. Similary, $\mathrm{ext}: \CT({\bf{X}}) \rightarrow \L^2({\bf{X}})$ will denote the natural extension operator. The $\H_\chi$ are closed subspaces, and the corresponding projection operators are  given by
\begin{equation*}
P_\chi= {d_\chi} \int_{G}  \overline{ \chi(k)} T(k) dk ,
\end{equation*}
 where $d_\chi$ is the dimension of the  irreducible representation corresponding to the character $\chi$, and $dk$ denotes Haar measure on $G$. If $G$ is just finite, $dk$ is the  counting measure, and one simply has
 \bqn
 P_\chi=\frac{d_\chi} {|G|} \sum_{k\in G} \overline{\chi(k)} T(k).
 \eqn
  Since $T(k)$ is unitary, one computes for $u,v \in \L^2(\R^n)$
 \begin{align*}
(u,P_\chi v) &=  {d_\chi} \int _{G}  \chi(k) ( u, T(k)v) dk = {d_\chi}  \int_{G} \overline{\chi (k^{-1})} (T(k^{-1})u,v) dk  =(P_\chi u, v),
\end{align*}
where we made use of $\overline{\chi(g)} =\chi (g^{-1})$. Hence $P_\chi$ is self-adjoint. Let now $A$ be the Friedrichs extension of the lower semi-bounded operator
\bqn
\mathrm{res} \circ A_0 \circ \mathrm{ext}: \CT({\bf{X}}) \longrightarrow \L^2({\bf{X}}).
\eqn
$A$ is a self-adjoint operator in $\L^2({\bf{X}})$, and is itself lower semi-bounded. Its spectrum is real, and consists of the point spectrum and the continuous spectrum. Recall that, in general, a symmetric operator $S$ in a separable Hilbert space is called lower semi-bounded, if there exists a real number $c$ such that
\bqn
( S u,u ) \geq  c \norm{u}^2 \qquad \text{for all } u \in \mathcal{D}(S),
\eqn
where $\mathcal{D}(S)$ denotes the domain of $S$. Now, if $V$ is a subspace contained in $ \mathcal{D}(S)$,  the quantity 
\bqn
\label{eq:15}
\Ncal(S,V)=\sup_{L\subset V}\{\dim L: ( S \,u,\,u) <0 \quad \forall \,\, 0\not=u \in \L\},
\eqn
can be used to give a qualitative description of the spectrum of $S$. More precisely, 
one has the following classical variational result of Glazman.
\begin{lemma}
Let $S$ be a self-adjoint, lower semi-bounded operator in a separable Hilbert space,
and define $N(\lambda,S)$ to be equal to the number of eigenvalues of $S$, counting multiplicities, less or equal $\lambda$, if  $(-\infty,\lambda)$ contains no points of the essential spectrum, and equal to $\infty$, otherwise. Then
\bqn
N(\lambda,S)=\Ncal(S-\lambda \1, \mathcal{D}(S)).
\eqn
\end{lemma}
\begin{proof} See \cite{levendorskii}, Lemma A.1.
\end{proof}
In particular, the lemma above allows one to determine whether $S$ has essential spectrum or not, where the latter is given by the continuous spectrum and the eigenvalues of infinite multiplicity.
Let us now return to the situation above. Since $A$ commutes with the action of $G$ on $\L^2({\bf{X}})$, the eigenspaces of $A$ are unitary $G$-modules that decompose into irreducible subspaces. Let therefore $N_\chi(\lambda)$ be equal to the number of eigenvalues of $A$, counting multiplicities, less or equal $\lambda$ and with eigenfunctions in $\mathrm{res }\, \H_\chi$,   if  $(-\infty,\lambda)$ contains no points of the essential spectrum, and equal to $\infty$, otherwise. One has then the following
\begin{lemma}
\label{lemma:8}
$ N_\chi(\lambda)=\Ncal(A_0-\lambda \1, \H_\chi \cap \CT({\bf{X}}))$.
\end{lemma}
\begin{proof}
Let $A_\chi$ be the Friedrichs extension of $\mathrm{res} \circ A_0 \circ \mathrm{ext}: \CT({\bf{X}}) \cap \H_\chi \longrightarrow \mathrm{res}\,  \H_\chi$. Then $N_\chi(\lambda)= N(\lambda, A_\chi)$, and the assertion follows with \cite{levendorskii}, Lemma A.2.
\end{proof}
 In order to estimate $\Ncal(A_0-\lambda \1, \H_\chi \cap \CT({\bf{X}}))$, we will apply the method of approximate spectral projection operators. It was first introduced by Tulovskii and Shubin,  and later developed and generalized by Feigin and Levendorskii, and we will mainly follow \cite{levendorskii} in our construction. Thus, let us  consider  on $\R^{2n}$ the metric
  \bq
 g_{x,\xi}(y,\eta)=|y|^2+ h(x,\xi)^2|\eta|^2 , \qquad h(x,\xi) = (1+|x|^2 +|\xi|^2 )^{-1/2}.
 \eq
which is clearly of the form  \eqref{eq:6}. Our symbol classes will be mainly of  the form   $S(h^{-2\delta}g, p)=\Gamma_{1-\delta,\delta}(p , \R^{2n})$ where $p$ is a $\sigma,h^{-2\delta}g$-temperate function, and $0 \leq \delta < 1/2$. In this case,
\bqn
h_\sigma^2(x,\xi)=(1 +|x|^2+|\xi|^2)^{2\delta-1},
\eqn
by equation \eqref{eq:2a}, 
which amounts to  $h_\sigma=h^{1-2\delta}$. 
Also note that $u \in S(h^{-2\delta} g, p)$ implies $\gd^\alpha_\xi \gd^\beta_x u \in S(h^{-2\delta}g, h^{(1-\delta)|\alpha| -\delta|\beta|}p)$. In particular, $S(h^{-2\delta} g, h^{-l})=\Gamma^l_{1-\delta,\delta}(\R^{2n})$, where $l\in \R$.   The symbols and functions used will also depend on the spectral parameter $\lambda$. Nevertheless, their membership to specific symbol classes will be uniform in $\lambda$, which  means that the values of their seminorms in the corresponding symbol classes will be bounded by some constant independent of $\lambda$. 
 Now, if $a$ denotes the left symbol of the classical pseudodifferential operator $A_0$, clearly $a \in S(g,h^{-2m}, K\times \R^n)$ for any compact set $K\subset \R^n$, so that $\sigma^l(A_0-\lambda\1) \in S(g,\tilde q_\lambda^2, K\times \R^n)$ uniformly in $\lambda \geq 1$, where
\bq
\label{eq:13}
\tilde q_\lambda^2(x,\xi)=h^{-2m} (x,\xi) +\lambda
\eq
is a $\sigma,g$-temperate function. But for $u \in \CT({\bf{X}})$, the quadratic form $((A_0-\lambda \1)u,u)$ entering in the definition of $\Ncal(A_0-\lambda \1, \H_\chi \cap \CT({\bf{X}}))$ depends only on values of $\sigma^l(A_0-\lambda \1)$ on ${\bf{X}} \times \R^n$. By changing  the latter symbol outside ${\bf{X}} \times \R^n$ we can achieve  that $\sigma^l(A_0-\lambda\1) \in S(g, \tilde q_\lambda^2)$ uniformly in $\lambda\geq 1$. In view of Corollary \ref{cor:1}  we can therefore assume that $A_0-\lambda \1$ can be represented as a pseudodifferential operator with Weyl symbol  $\tilde a_\lambda =\sigma^w(A_0-\lambda\1) \in S(g, \tilde q_\lambda^2)$. In particular, we may take  $\sigma^w(A_0) \in S(g,h^{-2m})$.  But by equation \eqref{eq:8a} and Lemma 13.1 in  \cite{levendorskii}  we  even have 
\bqn
a_{2m}(x,\xi) \geq c \qquad \text{ for all } (x,\xi) \in {\bf{X}} \times S^{n-1} \, \text{and some constant } c>0.
\eqn
Since $a -a_{2m} \in S(g,h^{-2m+1}, K \times \R^n)$, we can therefore  assume that   $A_0 \in \mathcal{LI}^+(g,h^{-2m})$, obtaining
\begin{lemma}
\label{lemma:11}
Let $A_0$ be a classical pseudodifferential operator satisfying \eqref{eq:8a}. Then $A_0$ and $A_0-\lambda \1$ can be represented as  pseudodifferential operators with Weyl symbols $\sigma^w(A_0)\in SI^+(g,h^{-2m})$ and  $\tilde a _\lambda \in {SI}^+(g,\tilde q_\lambda^2)$, respectively.
\end{lemma}
Note that if $\sigma^w(A_0)$, and consequently  also $\tilde a_\lambda$, are  $G$-invariant in the sense that 
\bqn
\sigma^w(A_0)( \sigma_g(x,\xi))=\sigma^w(A_0)(x,\xi), \qquad \tilde a_{\lambda} (\sigma_g(x,\xi))=\tilde a_{\lambda}(x,\xi),
\eqn
where  $\sigma_g$ is the symplectic transformation given by  $\sigma_g(x,\xi)=(\kappa_g(x), ^t \kappa _g '(x)^{-1}( \xi))=(\kappa_g(x),\kappa_g(\xi))$, and  $\kappa_g(x)=gx $ denotes the action of $g$,
the operators $A_0$ and $A_0 -\lambda\1$ will commute with the action of $G$ by Corollary \ref{cor:1}. We can therefore formulate the assumption about the $G$-invariance of $A_0$ also  in terms of its Weyl symbol, and shall henceforth assume that the Weyl symbol and the principal symbol $a_{2m}$ of $A_0$ are invariant under $\sigma_g$ for all $g\in G$.  In order to define the approximate spectral projection operators, we will introduce now the relevant symbols. Having in mind Lemma \ref{lemma:5}, let $a_\lambda \in S(g,1)$, and $d \in S(g,d)$ be $G$-invariant symbols which, on ${\bf{X}}_\rho \times \mklm{\xi: |\xi| >1}$, ${\bf{X}}_\rho = \mklm{ x: \dist(x, {\bf{X}}) <\rho}$, are given by
 \begin{align*}
  a_\lambda(x,\xi)&=\frac 1 {1+\lambda|\xi|^{-2m}} \Big ( 1 -\frac \lambda {a_{2m}(x,\xi)} \Big ),\\    d(x,\xi)&=|\xi|^{-1},
 \end{align*}
 where $\rho >0$ is some fixed constant, and in addition  assume that $d$ is positive and that $d(x,\xi)\to 0$ as $|x|\to \infty$. 
We  need to define smooth approximations to the Heaviside function, and to certain characteristic functions on  ${\bf{X}}$.
Thus, let $\tilde \chi$ be a smooth function on the real line satisfying 
$0 \leq \tilde \chi \leq 1$, and 
\bqn
\tilde \chi(s)=\left \{ \begin{array}{c} 1 \quad \text{for } s<0, \\ 0 \quad \text{for } s>1.
\end{array}\right.
\eqn
  Let $C_0>0$ and $\delta \in (1/4,1/2)$ be constants, and put $\omega=1/2 -\delta$. We then define the $G$-invariant function
\bq
\label{eq:16a}
\chi_\lambda=\tilde \chi\circ ((a_\lambda+4 h^{\delta-\omega}+8 C_0 d)\,  h^{-\delta}),
\eq
where $0 < \delta -\omega< 1/2$.
\begin{lemma}
\label{lemma:10a}
 $\chi_\lambda \in S(h^{-2\delta}g,1)=\Gamma^0_{1-\delta,\delta}(\R^{2n})$ uniformly in $\lambda$.
\end{lemma} 
\begin{proof}
We first note that 
\bqn
h^{-\delta} \in S(g,h^{-\delta}), \qquad (a_\lambda+4 h^{\delta-\omega}+ 8 C_0 d) \in S(g,1), 
\eqn
since  $d\in S(g,d) \subset S(g,1)$, and $h^{\delta-\omega} \in S(g,h^{\delta-\omega}) \subset S(g,1)$. Now, each of the derivatives of $\chi_\lambda$ with respect to $x$ and $\xi$ can be estimated by a sum of derivatives of $a_\lambda+4 h^{\delta-\omega}+8 C_0 d)\, h^{-\delta}$. But because of $\gd^\alpha_\xi\gd^\beta_x(a_\lambda+4h^{\delta-\omega}+8 C_0 d) \in S(g,h^{|\alpha|})$, $\gd^\alpha_\xi\gd^\beta_x h^{-\delta} \in S(g, h^{-\delta+|\alpha|})$, we obtain 
\bqn
|\gd^\alpha_\xi\gd^\beta_x \chi_\lambda(x,\xi)| \leq C_{\alpha,\beta} h^{(1-\delta)|\alpha|}= C_{\alpha,\beta} (1+|x|^2+|\xi|^2)^{-(1-\delta)|\alpha|/2},
\eqn
where $C_{\alpha,\beta}$ is independent of $\lambda$. We therefore obtain $\chi_\lambda \in \Gamma^0_{1-\delta,0}(\R^{2n}) \subset \Gamma^0_{1-\delta,\delta}(\R^{2n})$ uniformly in $\lambda$, and the assertion follows.
\end{proof}
Next, let $U$ be a subset in ${\R^{2n}}$, $c >0$, and put
\bqn 
U(c,g)= \mklm{(x,\xi)\in \R^{2n}: \exists   (y,\eta) \in U: g_{(x,\xi)} (x-y,\xi-\eta) < c};
\eqn 
according to Levendorskii \cite{levendorskii}, Corollary 1.2, there exists a smoothened characteristic  function $\psi_c\in S(g,1)$  belonging to the set $U$ and the parameter $c$, such that  $\supp \psi_c \subset U(2c,g)$, and ${\psi_c}_{|U(c,g)}=1$. 
 Let now 
 \bq
 \label{eq:19}
 \M_\lambda=\mklm{(x,\xi)\in \R^{2n}: a_\lambda< 4 h^{\delta-\omega} + 8 C_0 d}.
 \eq
 Both $\M_\lambda$ and $\gd{\bf{X}} \times \R^n$ are invariant under $\sigma_k$ for all $k \in G$, as well as  $(\gd {\bf{X}}  \times \R^n)(c,h^{-2\delta} g)$, and $\M_\lambda(c,h^{-2\delta}g)$, due to the invariance of $a_{2m}(x,\xi)$, and the considered metrics and symbols. Now, let $\tilde \eta_c, \psi_{\lambda,c}\in S(h^{-2\delta} g,1)$ be smoothened characteristic functions corresponding to the parameter $c$, and the sets $\gd {\bf{X}}  \times \R^n$ and $\M_\lambda$, respectively. According to  Lemma \ref{lemma:5}, we can assume that they are invariant under $\sigma_k$ for all $k \in G$; otherwise consider $ \int _{G} \tilde \eta_c \circ \,\sigma_k \,dk, \, \int_{ G} \psi_{\lambda,c} \circ \,\sigma_k \,dk $, respectively.
 We then define the functions
 \begin{align}
\eta_{\lambda,-c}(x,\xi)&= \left \{ \begin{array}{cc} 0, & x \notin {\bf{X}}, \\ ( 1-\tilde \eta_c(x,\xi)) \psi_{\lambda,1/c}(x,\xi), & x \in {\bf{X}}, \end{array} \right.\\
\eta_{c}(x,\xi)&= \left \{ \begin{array}{cc} \tilde \eta_c(x,\xi),  & x \notin {\bf{X}}, \\  1, & x \in {\bf{X}}. \end{array} \right.
\end{align}
Only  the support of $\psi_{\lambda,c}$ depends on $\lambda$, but not its growth properties, so that  $\eta_c,\eta_{\lambda,-c} \in S(h^{-2\delta}g,1)$ uniformly in $\lambda$.
Furthermore, since $\tilde \eta_{2c}=1$ on $\supp \tilde \eta_c$, and $\psi_{\lambda,1/c}=1 $ on $\supp \psi _{\lambda,1/2c}$, on has $\eta_{\lambda,-c}=1$ on $\supp \eta_{\lambda,-2c}$, which implies $\eta_{\lambda,-2c}\eta_{\lambda,-c} = \eta_{\lambda,-2c}.$ Similarly,  one verifies $\eta_c \eta_{2c}=\eta_c$. We are now ready to define the  approximate spectral projection operators.
  \begin{definition}
The approximate spectral projection operators of the first kind are defined by 
\bqn
\tilde \E_\lambda= \Op^w (\eta_{\lambda,-2})\, \Op^w(\chi_\lambda)\, \Op^w( \eta_{\lambda,-2}), 
 \eqn
while the approximate spectral projection operators of the second kind are
\bqn
\E_\lambda=\tilde \E_\lambda^2(3 -2\tilde \E_\lambda).
\eqn
\end{definition}
\begin{remark} $\tilde \E_\lambda$ is a smooth approximation to the spectral projection operator $E_\lambda$ of $A$ using  Weyl calculus, while $\E_\lambda$ is an approximation to $E_\lambda^2(3-2E_\lambda)=E_\lambda$. Note that, since $\eta_{\lambda,-2}$ and $\chi_\lambda$ are $G$-invariant, Corollary \ref{cor:1} implies that the operators $\Op^w (\eta_{\lambda,-2})$, $\Op^w(\chi_\lambda)$, and consequently also $\tilde \E_\lambda$ and $\E_\lambda$,  commute with the action $T(g)$ of $G$.  The choice of $\E_\lambda$ was originally  due to the fact that its trace class norm can be estimated from above by the operator norm of $3-2\tilde \E_\lambda$, and the Hilbert-Schmidt-norm of $\tilde \E_\lambda$, which are easier to handle. This construction was first used by Feigin \cite{feigin}.
\end{remark} 
Both $\tilde \E_\lambda$ and $\E_\lambda$ are integral operators with kernels in  $\S(\R^{2n})$. Indeed, the asymptotic expansion \eqref{eq:A}, together with Proposition  \ref{prop:A}, imply that the Weyl symbols of $\tilde \E_\lambda$ and $ \E_\lambda$ can be written in the form $a +r$, where $a$ has compact support, and $r \in S^{-\infty}(h^{-2\delta}g,1)$, because $\chi_\lambda$ has compact support in $\xi$, and $\eta_{\lambda,-2}$ has $x$-support in ${\bf{X}}$. Thus, $\sigma^w(\tilde \E_\lambda)$ and $\sigma^w(\E_\lambda)$  are rapidly decreasing Schwartz functions, and the same holds for the corresponding $\tau$-symbols. By Lemma 7.2 in \cite{hoermander79}, this also implies that $\tilde \E_\lambda$ and $\E_\lambda$ are of trace class and, in particular, compact operators in $\L^2(\R^n)$.  In addition, by Theorem \ref{thm:MT}, and the asymptotic expansion \eqref{eq:12bis},  one has $\sigma ^\tau(\tilde \E_\lambda)$, $\sigma^\tau(\E_\lambda)\in S(h^{-2\delta} g,1)$ uniformly in $\lambda$.
 On the other hand, the functions $\eta_{\lambda,-2}$ and $\chi_\lambda$ are real valued, which by general Weyl calculus implies that $ \Op^w (\eta_{\lambda,-2})$, $\Op^w(\chi_\lambda) $, and consequently also  $\tilde \E_\lambda$, and $\E_\lambda$,  are self-adjoint operators in  $\L^2(\R^n)$. By construction, $\E_\lambda$  commutes with the projection $P_\chi$,  so that  $P_\chi \E_\lambda=\E_\lambda P_\chi$ is a self-adjoint operator of trace class as well.
 Although the decay properties of $\sigma^\tau(\E_\lambda)$ are independent of $\lambda$, its support does depend on $\lambda$, which will result in  estimates for the trace of $P_\chi\E_\lambda$ in terms of  $\lambda$ that will be used in order to prove Theorem \ref{thm:1}. In particular, the estimate for the remainder term in Theorem \ref{thm:1} is determined by   the particular choice of the range $(1/4,/1/2)$ for the parameter $\delta$, which  guarantees that $1-\delta >\delta$.  By the general theory of compact, self-adoint operators, zero is the only accumulation point of the point spectra of $\tilde \E_\lambda$ and $\E_\lambda$, as well as the only point that could possibly belong to the continuous spectrum.
  The following proposition and its corollary give uniform bounds for the number of eigenvalues away from zero. They are based on certain $\L^2$-estimates for pseudodifferential operators.
\begin{proposition}
\label{prop:3}
The number of eigenvalues of $\tilde \E_\lambda$ lying outside  the interval $[-\frac 1 4, \frac 5 4]$ is bounded by some constant independent of $\lambda$.
\end{proposition}
\begin{proof}
Since $\chi_\lambda, \eta_{\lambda,-c} \in S(h^{-2\delta}g,1)$, Theorem \ref{thm:MT} yields $\sigma^w(\tilde \E_\lambda)\in  S(h^{-2\delta}g,1)$ uniformly in $\lambda$. Furthermore, taking into account the asymptotic expansion \eqref{eq:A}, we have 
\bqn
\sigma^w(\tilde \E_\lambda) =\eta^2_{\lambda,-2} \chi_\lambda +r_\lambda, 
\eqn
where $r_\lambda \in S(h^{-2\delta}g,  h^{1 - 2\delta})$. Now, since $0 \leq \chi_\lambda,\eta^2_{\lambda,-2}\leq 1$, for each $\epsilon >0$  there exists a constant $c>0$ such that $\epsilon +\eta^2_{\lambda,-2}\chi_\lambda \geq c$ and $(1+\epsilon) -\eta^2_{\lambda,-2}\chi_\lambda \geq c$. Consequently, the symbols of $\epsilon \1 +\tilde \E_\lambda$ and $(1+\epsilon)\1-\tilde \E_\lambda$ admit a representation of the form $a_1+a_2$, where $a_1 \geq c, \, a_2 \in S(h^{-2\delta}g, h^{1 - 2\delta})$; thus
\bqn
\epsilon \1+\tilde \E_\lambda \in \mathcal{LI}^+(h^{-2\delta}g,1), \qquad (1+\epsilon)\1 -\tilde \E_\lambda \in  \mathcal{LI}^+(h^{-2\delta}g,1)
\eqn
uniformly in $\lambda$.
According to   Lemma \ref{lemma:SI2}, this implies that for each $\lambda$ there  exist two operators $T_1,T_2$ such that $\epsilon \1 + \tilde \E_\lambda \geq T_1$ and $(1+\epsilon) \1-\tilde \E_\lambda \geq T_2$, and $T_i  \in \Lcal^{-\infty}(g,1)$ uniformly in $\lambda$. Therefore, by   Lemma \ref{lemma:3}, there exist  two subspaces $L_i\subset \L^2(\R^n)$ of finite codimension such that $\norm{T_i u}_{\L^2} \leq \epsilon \norm{u}_{\L^2}$ for $u \in L_i$ and all $\lambda$, which implies, via Cauchy-Schwartz, that $-\epsilon \norm{u}^2_{\L^2} \leq (T_iu,u) \leq \epsilon \norm{u}^2_{\L^2}$ on $L_i$. Putting everything together we arrive at the $\L^2$-estimates
\begin{align*}
(\tilde \E_\lambda u,u) &\geq ( (T_1-\epsilon \1)u,u) \geq -2\epsilon \norm{u}^2_{\L^2}, \\
(\tilde \E_\lambda u,u) &\leq ( ((1+\epsilon) \1-T_2)u,u) \leq (1+2\epsilon) \norm{u}^2_{\L^2},
\end{align*}
where $u \in L_1\cap L_2$, and taking $\epsilon=\frac 1 8 $ yields the desired result, since $\mathrm{codim} \,L_1 \cap L_2<\infty$.
\end{proof}
\begin{corollary}
\label{cor:2}
The number of eigenvalues of $\E_\lambda$ lying outside the interval $[0,1]$ is bounded by some constant independent of $\lambda$.
\end{corollary}
\begin{proof}
If $\tilde \nu_i$ denote the eigenvalues of $\tilde \E_\lambda$, then the eigenvalues of $\E_\lambda$ are given by $\nu_i=\tilde \nu_i^2(3-2\tilde \nu_i)$.
\end{proof}
Let now $N^{\E_\lambda}_\chi$ denote the number of eigenvalues of $\E_\lambda$ which are $\geq 1/2$, and whose eigenfunctions are contained in the $\chi$-isotypic component $\H_\chi$ of $\L^2(\R^n)$.  Since zero is the only accumulation point of the point spectrum of $\E_\lambda$,  $N^{\E_\lambda}_\chi$ is clearly finite. The next lemma will show that  it can be estimated by  the trace of the operator $P_\chi \E_\lambda$, and its square, so that it is natural to expect that it will  provide  a good approximation for $N_\chi(\lambda)=\tr P_\chi E_\lambda =\Ncal(A_0-\lambda \1, \H_\chi \cap \CT({\bf{X}}))$. 
\begin{lemma}\label{lemma:12a}
There exist constants $c_1,c_2>0$ independent of $\lambda$ such that
\bq
\label{eq:12a}
2 \tr (P_\chi \E_\lambda)^2 -\tr P_\chi \E_\lambda -c_1 \leq N^{\E_\lambda}_\chi \leq 3 \tr P_\chi \E_\lambda-2 \tr (P_\chi \E_\lambda)^2 +c_2. 
\eq
\end{lemma}
\begin{proof}
 Since $\E_\lambda\in \Lcal(h^{-2\delta}g,1)$, Theorem  \ref{thm:L2} implies that $\E_\lambda$ is $\L^2$-continuous. Moreover, by Remark \ref{rem:L2}, there is a constant C independent of $\lambda$ such that $\norm {\E_\lambda}_{\L^2}\leq C$; hence all eigenvalues of the operators $\E_\lambda$ are bounded by $C$. Let now  $\nu_{i,\chi}$ denote the eigenvalues of $\E_\lambda$ with eigenfunctions in $\H_\chi$. Taking into account Corollary \ref{cor:2} and the previous remark, we obtain the estimate
 \begin{align*}
 N^{\E_\lambda}_\chi &\leq \sum_{\nu_{i,\chi}\geq 1/2} \nu_{i,\chi}+\sum_{1/2 \leq \nu_{i,\chi}\leq 1} (1-\nu_{i,\chi}) +c_1
 \leq \sum_{\nu_{i,\chi}\geq 1/2} \nu_{i,\chi}+2\sum_{1/2 \leq \nu_{i,\chi}\leq 1} \nu_{i,\chi} (1-\nu_{i,\chi}) +c_1,
 \end{align*}
where $c_1>0$, like all other constants $c_i>0$ occurring in this proof, can be chosen independent of $\lambda$. Consequently, the right hand side can be estimated from above by 
 $ 3 \tr P_\chi \E_\lambda  -2 \tr P_\chi \E_\lambda \cdot P_\chi \E_\lambda +c_2 $.
 In the same way one computes
 \begin{align*}
 N^{\E_\lambda}_\chi&= \sum_{\nu_{i,\chi}\geq 1/2} \nu_{i,\chi}+\sum_{ \nu_{i,\chi}\geq 1/2 } (1-\nu_{i,\chi})\geq \sum _i \nu_{i,\chi} - \sum _{0 \leq \nu_{i,\chi} \leq 1/2}\nu_{i,\chi} -c_3\\
 &\geq  \sum _i \nu_{i,\chi} - 2 \sum _{0 \leq \nu_{i,\chi} \leq 1/2}\nu_{i,\chi}(1- \nu_{i,\chi}) -c_3 \geq  \sum _i \nu_{i,\chi} - 2 \sum _{i}\nu_{i,\chi}(1- \nu_{i,\chi}) -c_4, 
 \end{align*}
 where the right hand side can be estimated from below by $2 \tr P_\chi \E_\lambda \cdot P_\chi \E_\lambda -\tr P_\chi \E_\lambda  -c_4$. This completes the proof of \eqref{eq:38}.
\end{proof}
As the next section will show, $N^{\E_\lambda}_\chi$ will provide us with a lower bound for the spectral counting function $N_\chi(\lambda)$. Nevertheless, in order to obtain an upper bound as well, it will be necessary to introduce new approximations to the spectral projection operators. Namely, let 
\bqn
\chi_\lambda ^+=\tilde \chi(a_\lambda^+ h^{-\delta}), \qquad a_\lambda^+= a_\lambda -4 h^{\delta-\omega}-8C_0 d,
\eqn
where $\tilde \chi$ is defined as in \eqref{eq:16a}. As in Lemma \ref{lemma:10a}, one verifies that $\chi_\lambda ^+ \in S(h^{-2\delta}g,1)$ uniformly in $\lambda$.
\begin{definition}
The approximate spectral projection operators of the third kind are 
\bqn
\tilde \F_\lambda= \Op^w (\eta_{2}^2 \chi_\lambda^+),
 \eqn
while the approximate spectral projection operators of the fourth kind are
\bqn
\F_\lambda=\tilde \F_\lambda^2(3-2\tilde \F_\lambda).
\eqn
\end{definition}
Like the projection operators of the first and second kind, $\tilde \F_\lambda$ and $\F_\lambda$ are self-adjoint operators in $\L^2(\R^n)$ with kernels in $\S(\R^{2n})$, and therefore of trace class. Since $\F_\lambda$ commutes with $T(k)$, $P_\chi\F_\lambda$ is a self-adjoint operator of trace class, too.
Let  $M^{\F_\lambda}_\chi$ denote the number of eigenvalues of $\F_\lambda$ which are $\geq 1/2$, and whose eigenfunctions are contained in the $\chi$-isotypic compoment $\H_\chi$. Since Proposition \ref{prop:3} and Corollary \ref{cor:2} hold for $\tilde \F_\lambda$ and $\F_\lambda$ as well, we obtain 
   \begin{lemma}\label{lemma:12f}
There exist constants $c_1,c_2>0$ independent of $\lambda$ such that
\bq
\label{eq:37}
2 \tr (P_\chi \F_\lambda)^2 -\tr P_\chi \F_\lambda -c_1 \leq M^{\F_\lambda}_\chi \leq 3 \tr P_\chi \F_\lambda-2 \tr (P_\chi \F_\lambda)^2 +c_2. 
\eq
\end{lemma}
\begin{proof}
The proof is a verbatim repetition of the proof of Lemma \ref{lemma:12a} with $\E_\lambda$ replaced by $\F_\lambda$.
\end{proof}

\section{Estimates from below for $N_\chi(\lambda)$}

In this section, we  shall estimate the spectral counting function $N_\chi(\lambda)=\mathcal{N}(A_0-\lambda \1,\H_\chi \cap \CT({\bf{X}}))$ from below. More precisely, by adapting techniques developed in \cite{levendorskii} to our situation, we will show the following 
 \begin{theorem}
 \label{thm:4}
Let $N^{\E_\lambda}_\chi$  be the number of eigenvalues of $\E_\lambda$ which are $\geq 1/2$, and whose eigenfunctions are contained in the $\chi$-isotypic component $\H_\chi$.
Then there exists a constant $C>0$ independent of $\lambda$ such that
\bq
 \Ncal(A_0-\lambda \1, \H_\chi \cap \CT({\bf{X}}))\geq N^{\E_\lambda}_\chi - C.
\eq
 \end{theorem} 
As a first step towards the proof, let $\tilde q_\lambda $ be defined as in \eqref{eq:13}, and   $q_\lambda \in SI(g, \tilde q_\lambda^{-1})$ be a $G$-invariant symbol which, on ${\bf{X}}_\epsilon \times \mklm{\xi:|\xi|>1}$ is given by 
\bqn
q_\lambda(x,\xi)=\big ( a_{2m}(x,\xi) (1+|\xi|^{-2m} \lambda)\big )^{-1/2},
\eqn
and consider the $G$-invariant function  $\pi=(h^{\delta-\omega} + C_0 d)^{-1/2}\in SI(g,\pi)$, together with  the operators
 \bqn
 \Pi= \Op^w(\pi), \qquad Q_\lambda =  \Op^w(q_\lambda).
 \eqn
 Since $\pi\tilde q_\lambda^{-1}$ is bounded, $\Pi Q_\lambda$ is a  continuous operator in $\L^2(\R^n)$. The parametrices of $\Pi$ and $Q_\lambda$, which exist according to Lemma \ref{lemma:SI1}, will be denoted  by  $R_\Pi$ and $R_{Q_\lambda}$. Furthermore, an examination of the proof of Lemma \ref{lemma:SI1} shows that if $a\in SI(g,m)$ is $G$-invariant, then the Weyl symbol $b$ of the parametrix of $\Op^w(a)$ can be assumed to be $G$-invariant. Consequently, the parametrices  $R_\Pi$ and $R_{Q_\lambda}$ commute with the operators $T(k)$.
 \begin{lemma}
 Let $L_\chi^{\E_\lambda}=\mathrm{Span} \{ u \in \S(\R^n) \cap \H_\chi: \E_\lambda u = \nu u,  \nu \geq \frac 1 2\}$ and  $\tilde L_\chi^{\E_\lambda}= \Op^l(\eta_{\lambda,-1})\, Q_\lambda \,\Pi L_\chi^{\E_\lambda}$. Then
 \bq
 \dim \tilde L_\chi^{\E_\lambda} \geq \dim L_\chi^{\E_\lambda} - C =N^{\E_\lambda}_\chi - C
 \eq
 for some constant $C>0$ independent of $\lambda$.
 \end{lemma}
 \begin{proof}
 Let us first note that since  $\eta_{\lambda,-1}$ has support in ${\bf{X}} \times \R^n$, and $\Op^l(\eta_{\lambda,-1})\, Q_\lambda \,\Pi$ commutes with $P_\chi$, we have $\tilde L_\chi^{\E_\lambda} \subset \CT({\bf{X}}) \cap \H_\chi$. 
 Next, we will prove that 
 \bq
 \label{eq:42}
 R_\Pi R_{Q_\lambda} \, \Op^l(\eta_{\lambda,-1}) \, Q_\lambda \,\Pi \E_\lambda = \E_\lambda+T, 
 \eq
 where $T \in \mathcal L{^{-\infty}}(g,1)$. Indeed, the Weyl symbol of $ \Op^l(\eta_{\lambda,-1}) \, Q_\lambda \,\Pi \E_\lambda$ is given by a linear combination of products of derivatives of the Weyl symbols of $Q_\lambda$, $\Pi$, $\E_\lambda$, and $\Op^l(\eta_{\lambda,-1})$. By  the asymptotic expansion \eqref{eq:12bis},
   \begin{align*}
 \sigma^w(\Op^l(\eta_{\lambda,-1}))(x,\xi) &\sim  \sum_{\beta} \frac  1 {\beta !} \Big (\frac {-1} 2 \Big )^{|\beta|} \gd _\xi^{\beta}  D_x^\beta  \eta_{\lambda, -1}(x,\xi).
 \end{align*} 
Now, equation \eqref{eq:23b} implies that, up to terms of order $-\infty$, the support of $ \sigma^w(\E_\lambda)$ is contained in $\supp \eta_{\lambda,-2}$, and we  shall express this by writing $\supp_\infty \sigma^w(\E_\lambda) \subset \supp \eta_{\lambda,-2}$. For the same reason,  we must have $\supp_\infty \sigma^w(\Op^l(\eta_{\lambda,-1}) \, Q_\lambda \,\Pi \E_\lambda)\subset \supp \eta_{\lambda,-2}$. But $\eta_{\lambda,-1}=1$ on $\supp \eta_{\lambda,-2}$ implies that all terms in the expansion of $\sigma^w(\Op^l(\eta_{\lambda,-1}))$ vanish on $\supp \eta_{\lambda,-2}$, except for the zero order terms. Proposition \ref{prop:A} then yields 
 \bqn
\sigma^w(\Op^l(\eta_{\lambda,-1}))(x,\xi) =\eta_{\lambda,-1}(x,\xi) 
 \eqn 
on $ \supp \eta_{\lambda,-2}$, up to a term of order $-\infty$. On this  set, the Weyl symbol of $ \Op^l(\eta_{\lambda,-1}) \, Q_\lambda \,\Pi \E_\lambda$  therefore reduces to   $\eta_{\lambda,-1}=1$ times a linear combination of products of derivatives of the Weyl symbols of $Q_\lambda$, $\Pi$ and $\E_\lambda$ supported in $\supp \eta_{\lambda,-2}$, which corresponds to the Weyl symbol of $Q_\lambda \,\Pi \E_\lambda$, plus  an additional term of order $-\infty$. Thus,
\bq
\label{eq:43}
 \Op^l(\eta_{\lambda,-1}) \, Q_\lambda \,\Pi \E_\lambda= Q_\lambda \,\Pi \E_\lambda +\tilde T, \qquad \tilde T \in \Lcal^{-\infty}(g,\pi \tilde q_\lambda^{-1}),
\eq
and \eqref{eq:42} follows by taking into account the definition of the parametrix.
 Now, $\E_\lambda:L_\chi^{\E_\lambda}\rightarrow L_\chi^{\E_\lambda}$ is clearly surjective, and 
 \bqn
 \norm{\E_\lambda u} \geq \frac 1 2 \norm{u}, \qquad u \in L_\chi^{\E_\lambda},
 \eqn
 implies that $\E_\lambda$ is injective on $L_\chi^{\E_\lambda}$ as well. Equation \eqref{eq:42} therefore means that on $L_\chi^{\E_\lambda}$
 \bq
 \label{eq:44}
  R_\Pi R_{Q_\lambda} \, \Op^l(\eta_{\lambda,-1}) \, Q_\lambda \,\Pi= \1_{L_\chi^{\E_\lambda}} + T \E_\lambda^{-1}.
 \eq
 According to Lemma \ref{lemma:3}, there exists a subspace of finite codimension $M$ such that $\norm{Tu} \leq \norm{u}/8$ for all $ u \in M$ and all $\lambda$. This gives
 \bqn
 \norm{T\E_\lambda^{-1}u}\leq 2 \norm{Tu}\leq \frac 1 4 \norm{u} \qquad \text{for all  } u \in L_\chi^{\E_\lambda}\cap M.
 \eqn
 Let now $v,w \in L_\chi^{\E_\lambda} \cap M$, and assume that $\Op^l(\eta_{\lambda,-1}) \, Q_\lambda \,\Pi v = \Op^l(\eta_{\lambda,-1}) \, Q_\lambda \,\Pi w$. By \eqref{eq:44} we deduce
 $w+T\E_\lambda^{-1} w=v+T\E_\lambda^{-1}v$ and consequently $\norm{(\1+T\E_\lambda^{-1})(v-w)}=0$. But for $u\in M \cap L_\chi^{\E_\lambda}$ one computes
 \bqn 
 \norm{ ( \1 + T\E_\lambda^{-1}) u} \geq \norm {u}-\norm{T\E_\lambda^{-1}u } \geq \Big ( 1 - \frac 1 4 \Big ) \norm{u};
 \eqn
hence $\1+T\E_\lambda^{-1}$ is injective, and $v=w$. Thus we have shown that
\bq
\label{eq:45}
\Op^l(\eta_{\lambda,-1}) \, Q_\lambda \,\Pi: L_\chi^{\E_\lambda} \cap M \longrightarrow \tilde L_\chi^{\E_\lambda}
\eq 
is injective, and the assertion of the lemma follows with $C=\mathrm{codim} \,M<\infty $.
 \end{proof}
 Since $\tilde L_\chi^{\E_\lambda} \subset \CT ({\bf{X}}) \cap \H_\chi$, the next proposition will provide us with a suitable reference subspace in order  to prove Theorem \ref{thm:4}. Its dimension will  be estimated from below with the help of the preceding lemma. Note that the parametrices of $\Pi$ and $Q_\lambda$ were needed to show the injectivity of \eqref{eq:45}. 
 \begin{proposition}
 There exists a subspace $L\subset \tilde L_\chi^{\E_\lambda}$  such that $\dim L \geq \dim L_\chi^{\E_\lambda}-C$ for some constant $C>0$ independent of $\lambda$, and
 \bqn
 ((A_0-\lambda \1)u,u) _{\L^2} < 0 \quad \text{for all } 0 \not= u \in L.
 \eqn
  \end{proposition}
  Note that, by construction, $\tilde L_\chi^{\E_\lambda} \subset \CT({\bf{X}}) \cap \H_\chi$, while $L_\chi^{\E_\lambda} \not\subset \CT({\bf{X}})$. It is  this proposition that accomplishes the transition from $\R^n$ to ${\bf{X}}$, which, according to \eqref{eq:43}, is achieved by a perturbation of order $-\infty$.
  
 \begin{proof}
 Let $ v \in L_\chi^{\E_\lambda}$ and $w=\Op^l(\eta_{\lambda,-1}) \, Q_\lambda \,\Pi \E_\lambda v \in \tilde L_\chi^{\E_\lambda}$. Equation \eqref{eq:43} implies that
 \bqn 
 w = Q_\lambda \,\Pi  \E_\lambda v + \tilde T v, \qquad \tilde T \in \Lcal^{-\infty}(g,\pi \tilde q_\lambda^{-1}).
 \eqn
 Consequently, one computes
 \begin{align}
 \label{eq:47}
 \begin{split}
 ((A_0-\lambda \1)w,w)&=\big (\Pi^\ast \, Q_\lambda^\ast [ A_0-\lambda \1 + 4 R_{Q_\lambda}^\ast  \Op^w(h^{\delta - \omega}+ C_0 d)R_{Q_\lambda}] \, Q_\lambda \,\Pi \E_\lambda v, \E_\lambda v\big )\\
 &-4( \Pi^\ast \, Q_\lambda^\ast R_{Q_\lambda}^\ast  \Op^w (h^{\delta-\omega}+ C_0 d)R_{Q_\lambda} \, Q_\lambda \,\Pi \E_\lambda v, \E_\lambda v) +(Tv,v)
 \\&=:(D_1 \E_\lambda v, \E_\lambda v) - 4 (D_2 \E_\lambda v, \E_\lambda v) + ( T v,v),
  \end{split} 
  \end{align}
 where $ T $ is of order $-\infty$. Now, since $Q_\lambda R_{Q_\lambda} - \1 \in \Lcal^{-\infty}(g,1)$, we have
\bqn
D_2 = \1 + K_2, \qquad K_2 \in \Lcal(g,h); 
\eqn
 indeed, by definition, the Weyl symbol of $\Pi$ is equal to $\pi=(h^{\delta-\omega}+ C_0 d )^{-1/2}\in SI(g,\pi)$. Now, according to Lemma \ref{lemma:11},  $A_0-\lambda \1=\Op^w(\tilde a _\lambda)$, where $\tilde a_\lambda \in SI^+(g, \tilde q_\lambda^2)$. Thus, 
 \bq
 D_1 = \,\Pi^\ast [ Q_\lambda^\ast \, \Op^w(\tilde a_\lambda) \, Q_\lambda + 4 \, \Op^w (h^{\delta -\omega} +C_0 d) ] \,\Pi + K_1,
 \eq
 where $K_1 \in \Lcal^{-\infty}(g,1)$.  Furthermore, we can assume that $ q_\lambda \in S(g,\tilde q_\lambda^{-1})$ is such that $a_\lambda=(a_{2m}-\lambda)q^2_\lambda \in S(g,1)$, and using Theorem \ref{thm:MT}  one computes
 \bq
 \label{eq:41a}
  a_\lambda - \sigma^w(Q_\lambda^\ast \, \Op^w(\tilde a_\lambda) \, Q_\lambda) = a_\lambda - q_\lambda^2 \tilde a_\lambda + r =  q_\lambda^2(a_{2m}-\lambda-\tilde a_\lambda) + r \in S(g,d),
  \eq
  where $r \in S(g,h)$. But this implies $a_\lambda -\sigma^w(Q_\lambda^\ast \Op^w(\tilde a_\lambda) Q_\lambda) + 4 C_0d \geq c d$ for some sufficiently large $C_0$ and some $c>0$; hence
  \bq
  \label{eq:49}
  a_\lambda - \sigma^w(Q_\lambda^\ast \, \Op^w(\tilde a_\lambda) \, Q_\lambda) +4 C_0 d \in SI^+(g,d).
  \eq
 Using Lemma \ref{lemma:SI2}, we conclude from \eqref{eq:49} that there exists a $T_4 \in \Lcal^{-\infty}(g,d)$ such that
 \bq
 \label{eq:50}
 Q_\lambda^\ast \, \Op^w(\tilde a_\lambda) \, Q_\lambda \leq \Op^w(a_\lambda) + 4 C_0 \, \Op^w(d) + T_4.
  \eq
Together with $\norm {\E_\lambda v }^2 \geq \frac 1 4 \norm{v}^2$, equations \eqref{eq:47} - \eqref{eq:50} therefore yield the estimate
  \begin{align*}
  ((A_0-\lambda \1) w,w ) &= (Tv,v)- 4(K_2 \E_\lambda v, \E_\lambda v) - 4(\E_\lambda v,\E_\lambda v) + (K_1 \E_\lambda v,\E_\lambda v)\\
  &+ ( \Pi^\ast [Q_\lambda^\ast \, \Op^w (\tilde a_\lambda) \, Q_\lambda + 4C_0 \, \Op^w(d) + 4 \, \Op^w (h^{\delta-\omega})] \,\Pi \E_\lambda v, \E_\lambda v) \\
  &\leq ( \Pi^\ast [ \Op^w (a_\lambda) +8 C_0 \, \Op^w(d) + 4 \, \Op^w(h^{\delta-\omega})] \,\Pi \E_\lambda v, \E_\lambda v)-\norm{v}^2 + (K_3 v,v),
  \end{align*}
  where $K_3 \in S(h^{-2\delta}g,h)$. We therefore set
  \bqn
  a_\lambda^-:= a_\lambda + 8 C_0 d +4 h^{\delta-\omega} \in S(g,1),
  \eqn
  and obtain the estimate
  \bq
  \label{eq:51}
  ((A_0-\lambda \1) w,w ) \leq ( \Pi^\ast  \, \Op^w (a_\lambda^-)  \,\Pi \E_\lambda v, \E_\lambda v)-\norm{v}^2 + (K_3 v,v).
  \eq
  Thus, it remains to show that $\E_\lambda^\ast \,\Pi^\ast \, \Op^w (a_\lambda^-) \,\Pi \E_\lambda -\1+K_3$ is negative definite on some subspace of finite codimension. In order to do so, we will show that $\E_\lambda^\ast \,\Pi^\ast \, \Op^w (a_\lambda^-) \,\Pi \E_\lambda -\1+K_3 \leq -\1 + K_4$, where $K_4 \in \Lcal(h^{-2\delta} g, h^\omega)$ and  $\omega >0$. As it shall become apparent in the following discussion, the key to this is contained in  the fact that, although $a_\lambda^- \in S(g,1)$, there exists a $K_5 \in \Lcal(h^{-2\delta}g, h^\delta)$ such that $\Op^w(\chi_\lambda a_\lambda^- \chi_\lambda) \leq K_5$! Now,
  \begin{align*}
  \Pi \E_\lambda &= \Pi \tilde \E_\lambda \D_\lambda =\Pi  \, \Op^w(\eta_{\lambda,-2})\, \Op^w (\chi_\lambda)\, \Op^w(\eta_{\lambda,-2})  \D_\lambda \\
&= [\Pi \, \Op^w(\eta_{\lambda,-2}),\, \Op^w(\chi_\lambda)] \, \Op^w(\eta_{\lambda,-2})  \D_\lambda \\
  & +  \Op^w(\chi_\lambda) \,\Pi \, \Op^w(\eta_{\lambda,-2}) \, \Op^w(\eta_{\lambda,-2})) \D_\lambda =: W_1 +W_2,
  \end{align*}
  where we put $\D_\lambda= \tilde \E_\lambda ( 3 -2 \tilde \E_\lambda)$. Since $\Pi$ and $\tilde \E_\lambda$ are self-adjoint, we obtain
\begin{align}
\label{eq:48}
\E_\lambda &\,\Pi \, \Op^w (a_\lambda^-) \,\Pi \E_\lambda  =  W_2^\ast  \, \Op^w(a_\lambda^-) W_2 + R 
\end{align}  
  where $R=W_1^\ast \, \Op^w(a_\lambda^-)W_2 + W_2^\ast \, \Op^w (a_\lambda^-) W_1 + W_1^\ast \, \Op^w(a_\lambda^-) W_1$ is given by a sum of terms which contain either $[\Pi \, \Op^w(\eta_{\lambda,-2}), \, \Op^w(\chi_\lambda)]$, or its adjoint $[\Op^w(\chi_\lambda),\, \Op^w(\eta_{\lambda,-2})\,\Pi]$, as factors. Now, the crucial remark is  that 
  \bq
  \label{eq:44a}
   \supp_\infty \sigma ^w ([ \Pi \, \Op^w(\eta_{\lambda,-2}),\, \Op^w(\chi_\lambda)] ) \subset \supp_{\mathrm{diff}} \chi_\lambda \subset \mklm{(x,\xi): |a_\lambda^-(x,\xi)| \leq h^\delta(x,\xi)},
  \eq
  where $\supp_{\mathrm{diff}} \chi_\lambda=\mklm{ (x,\xi): \exists k>0: \chi_\lambda^{(k)}(x,\xi) \not = 0}$. 
  To see this, first note that by Theorem \ref{thm:MT} and Proposition \ref{prop:A}, we have the trivial inclusion $\supp_\infty \sigma ^w ([ \Pi \, \Op^w(\eta_{\lambda,-2}),\, \Op^w(\chi_\lambda)] )\subset \supp \chi_\lambda.$ But since the terms in the asymptotic expansion of the Weyl symbol of $[ \Pi \, \Op^w(\eta_{\lambda,-2}),\, \Op^w(\chi_\lambda)]$ are of order $\geq 1$, they vanish unless $(x,\xi) \in \supp _{\mathrm{diff}} \chi_\lambda$, and one obtains the first inclusion. The second inclusion follows by noting the implications
  \bqn
   \chi_\lambda^{(k)} =0 \,\forall\, k>0 \quad \Leftarrow \quad \chi_\lambda=0 \text{ or } \chi_\lambda=1 \quad \Leftarrow \quad  a_\lambda^- h^{-\delta} \geq 1 \text{ or } a_\lambda^- h^{-\delta} \leq 0.
  \eqn
  While computing the Weyl symbol of $R$, we can therefore replace $a_\lambda^-$ with 
  \bq
  \label{eq:52c}
  b_\lambda^-=a_\lambda^-  \theta_\lambda, \qquad \theta_\lambda =\theta \Big( \frac 1 2 a_\lambda^-h^{-\delta}\Big ),
  \eq
  where $\theta\in \CT(\R)$ is a real valued function taking values between $0$ and $1$, which is equal $1$ on $[-1,1]$, and which vanishes outside $[-2,2]$, so that $\theta_\lambda=1$ on $\mklm{(x,\xi): |a_\lambda^-(x,\xi)|\leq h^\delta(x,\xi) }$. Indeed, this replacement adds at most a term of order $-\infty$ to the Weyl symbol of $R$. Now, the adventage of performing this replacement resides in the fact that, on $\supp \theta_\lambda$, one has $|a_\lambda^- | \leq 4 h^\delta$, which together with
 \bqn
 |\gd_\xi^\alpha \gd _x^\beta a_\lambda^-(x,\xi)| \leq C(1+|x|^2+|\xi|^2) ^{\frac{-|\alpha|}2} \leq  C(1+|x|^2+|\xi|^2) ^{\frac{-\delta -(1-\delta)|\alpha|+\delta |\beta|}2}, \quad |\alpha| \geq 1,
 \eqn  
i.e.  $\nu_k(h^{-2\delta} g, h^\delta;a_\lambda^-) < \infty$, $k \geq 1$, yields $a_\lambda^- \in S (h^{-2\delta} g, h^\delta, \supp \theta_\lambda)$, in contraposition to $a_\lambda^- \in  S(g,1)$. Consequently,
 $ b_\lambda^- \in S (h^{-2\delta} g, h^\delta)$, and we obtain
  \bq
  \label{eq:52}
  R \in \Lcal(h^{-2\delta} g, h^\delta \pi^2) \subset \Lcal(h^{-2\delta}g, h^\omega),
  \eq
  since $W_1,W_2 \in \Lcal (h^{-2\delta}g,\pi)$, $ \D_\lambda \in \Lcal(h^{-2\delta}g,1)$, and $h^\delta\pi^2 =h^\delta(h^{\delta-\omega}+C_0 d)^{-1}= (h^{-\omega}+ C_0h^{-\delta} d)^{-1} \sim h^\omega$.
Equations  \eqref{eq:51}, \eqref{eq:48}, and \eqref{eq:52} therefore yield the estimate
 \bq
 \label{eq:53} ((A_0-\lambda \1) w,w) \leq (W_2^\ast \, \Op^w( a_\lambda^-) W_2 v,v) - \norm {v}^2 + (K_4v,v),
 \eq
where $K_4 =K_3 + R \in \Lcal (h^{-2\delta}g, h^\omega)$. To examine $W_2^\ast \Op^w(a_\lambda^-) W_2$ more closely, let us consider  the operator
\bqn
S= \Op^w (\chi_\lambda) \Op^w(a_\lambda^-)  \Op^w(\chi_\lambda) -\, \Op^w(\chi_\lambda a_\lambda^- \chi_\lambda).
\eqn
By the usual argument, the asymptotic expansion \eqref{eq:A} and  Proposition \ref{prop:A} yield $\supp_\infty \sigma^w(S) \subset \supp_{\mathrm{diff}} \chi_\lambda$. In the computation of the Weyl symbol of $S$ we can therefore again replace $a_\lambda^-$ with $b_\lambda^-$, getting at most an additional term of order $-\infty$. Since $ \Op^w(\chi_\lambda) \in \Lcal(h^{-2\delta}g,1)$ by Lemma \ref{lemma:10a}, we obtain
\bq
\label{eq:54}
S  \in \Lcal(h^{-2\delta}g, h^\delta).
\eq
Now, by construction, $ a_\lambda^-\chi_\lambda \leq h^\delta$, since $0 \leq \chi_\lambda \leq 1$ and $\chi_\lambda=0$ for $a_\lambda^-h^{-\delta} > 1$, so that one infers
\bqn
|\gd_\xi^\alpha\gd_x ^\beta ( \chi_\lambda a_\lambda^- \chi_\lambda)(x,\xi)| \leq C (1 + |x|^2+|\xi|^2) ^{(-\delta-(1-\delta)|\alpha|+\delta|\beta|)/2}
\eqn
for some constant $C>0$. 
 But this implies $\Op^w(\chi_\lambda a_\lambda^- \chi_\lambda) \in \Lcal(h^{-2\delta}g,h^\delta)$. Using \eqref{eq:53} and \eqref{eq:54}, we therefore get
\bqn
((A_0-\lambda \1)w,w) \leq ( W_3^\ast \, \Op^w(\chi_\lambda a_\lambda^-\chi_\lambda) W_3 v,v) - \norm{v}^2+ (K_5 v,v) = -\norm{v}^2+(K_6 v,v), 
\eqn
with
\begin{align*}
W_3&= \Pi \, \Op^w(\eta_{\lambda,-2})\, \Op^w(\eta_{\lambda,-2})  \D_\lambda \in \Lcal(h^{-2\delta} g,\pi),\\
K_5&=K_4 + W_3^\ast [\Op^w (\chi_\lambda) \Op^w(a^-_\lambda)\Op^w (\chi_\lambda) -  \Op ^w(\chi_\lambda a^-_\lambda \chi_\lambda)] W_3 \in \Lcal (h^{-2\delta } g, h^\omega), \\
K_6&=K_5 +W_3^\ast \Op^w(\chi_\lambda a_\lambda^- \chi_\lambda) W_3 \in \Lcal(h^{-2\delta} g, h^\omega).
\end{align*}
Since $h_\sigma= h^{1-2\delta}$, Lemma \ref{lemma:3} implies that  the operator $-\1+K_6$ is negative definite on a subspace $U \subset \L^2(\R^n)$ of finite codimension which does not depend on $\lambda$. Putting $L:=\Op^l(\eta_{\lambda,-1}) \, Q_\lambda \,\Pi \E_\lambda(U \cap L_\chi^{\E_\lambda}\cap M)\subset \tilde L_\chi^{\E_\lambda}$ with $M$ as in \eqref{eq:45}, we finally get
\bq
\label{eq:55}
((A_0-\lambda \1)w,w) < 0  \quad \forall \, 0 \not=w \, \in L,
\eq
where  $\dim U \cap L_\chi^{\E_\lambda} \cap M -\mathrm{codim} \,M  \leq  \dim M\cap \E_\lambda(U \cap L_\chi^{\E_\lambda} \cap M)   \leq \dim L $,  since $\E_\lambda$ is bijective on $L_\chi^{\E_\lambda}$, and $ \dim L_\chi^{\E_\lambda} \leq \dim U\cap L_\chi^{\E_\lambda} \cap M + \mathrm{codim} \, U\cap M $. The assertion of the proposition now follows.
\end{proof}
We can now prove Theorem \ref{thm:4}.
\begin{proof}[Proof of Theorem \ref{thm:4}]
Let $L\subset \tilde L_\chi^{\E_\lambda}\subset \CT({\bf{X}}) \cap \H_\chi$ be as in the  previous proposition. Then  \eqref{eq:55} holds, and $\mathcal{N}(A_0-\lambda \1,\H_\chi\cap \CT({\bf{X}})) \geq \dim L$. Furthermore, $\dim L \geq \dim L_\chi^{\E_\lambda} -C =N^{\E_\lambda}_\chi -C$, and the assertion of the theorem follows.
\end{proof}

\section{Estimates from above for $N_\chi(\lambda)$}

In this section, we will prove an estimate from above for $N_\chi(\lambda)=\mathcal{N}(A_0-\lambda \1,\H_\chi \cap \CT({\bf{X}}))$ in terms of 
the number $M_\chi^{\F_\lambda}$ of eigenvalues of $\F_\lambda$  which are $\geq 1/2$, and whose eigenfunctions are contained in the $\chi$-isotypic component $\H_\chi$.
 In order to do so,  we first prove the following
\begin{proposition}
\label{prop:8}
There exists a constant $C>0$ independent of $\lambda$ such that
\bqn
\Ncal(A_0-\lambda \1, \H_\chi \cap \CT({\bf{X}})) \leq \mathcal{N}(\Op^l(\eta_{1}) \,\Pi \, \Op^w(a_\lambda^+) \,\Pi \, \Op^r(\eta_{1})+\1,\H_\chi\cap \CT(\R^n)) +C.
\eqn
\end{proposition}
Note that this proposition accomplishes the transition from variational quantities related to $\R^n$ to quantities related to the bounded subdomain ${\bf{X}}$.  Now, the proof of Proposition \ref{prop:8} relies on the following
 \begin{lemma}
 \label{lemma:10}
 There exists a subspace $L \subset \CT({\bf{X}})$ of finite codimension in $\CT({\bf{X}})$ such that
 \bqn
 ( (A_0-\lambda \1) \,u, \,u) \geq ( [ \Op^l(\eta_{1}) \,\Pi \, \Op^w(a_\lambda^+ )\,\Pi \, \Op^r(\eta_{1})\, +\1]\, R_\Pi\, R_{Q_\lambda}\, u, \,R_\Pi \,R_{Q_\lambda}\, u)
 \eqn
 for all $0\not=u \in L$, and all $\lambda$.
 \end{lemma}
 \begin{proof}[Proof of  Proposition \ref{prop:8}]
 Let us assume Lemma \ref{lemma:10} for a moment, and introduce the notation
 \begin{align*}
 A_\lambda[u]=(A_0-\lambda \1)u,u), \quad B_\lambda[u]= ([\Op^l(\eta_{1}) \,\Pi \, \Op^w(a_\lambda^+ )\,\Pi \, \Op^r(\eta_{1})\, +\1]u,u).
 \end{align*}
According to that lemma,  there exists    a subspace $L$ in $\CT({\bf{X}})$ 
of finite codimension such that
\bqn
A_\lambda[u] \geq B_\lambda[R_\Pi R_{Q_\lambda}u], \quad 0 \not= u \in L,
\eqn
for all $\lambda$. Let now $m $ be a positive $\sigma, g$ -temperate function such that $1/m$ is bounded. Following \cite{levendorskii}, we introduce the weight spaces of Sobolev type
\bqn
\H(g,m) =\mathrm{span} \mklm{ Tw:w \in \L^2(\R^n), \, T \in \Lcal(g,1/m) } \subset \L^2(\R^n),
\eqn
and endow them with the strongest topology in which each of the operators $T:\L^2(\R^n)\rightarrow \H(g,m)$, $ T \in \Lcal(g,1/m)$, is continuous. It can then be shown that there exists an operator $\Lambda_m \in \Lcal(g,m)$ such that $\Lambda_m:\H(g,m) \rightarrow \L^2(\R^n)$ is a topological isomorphism. In particular, $\H(g,m)$ becomes a Hilbert space with the norm $\norm{u}_m=\norm{\Lambda_mu}_{\L^2}$. Furthermore, we have the continuous embedding $\S(\R^n) \subset \H(g,m)$, and if $m_1$ is a bounded, $\sigma,g$-tempered function, and $A \in \L(g,mm_1)$, then $A:\H(g,m) \rightarrow \H(g,m_1^{-1})$ defines  a continuous map. Now, by Theorem \ref{thm:3} and the asymptotic expansion \eqref{eq:12}, $R_\Pi R_{Q_\lambda} \in \mathcal{LI}(g, \pi^{-1} \tilde q_\lambda)$, so that by Lemma \ref{lemma:SI1} the operator $\Lambda_\pi R_\Pi R_{Q_\lambda}  \Lambda^{-1}_{\tilde q_\lambda}\in \mathcal{LI}(g,1)$ has a parametrix $Z \in \mathcal{LI}(g,1)$ satisfying $Z \, \Lambda_\pi R_\Pi R_{Q_\lambda}  \Lambda^{-1}_{\tilde q_\lambda}= \1 +K$, where $K \in \Lcal^{-\infty}(g,1)$. Since by Lemma \ref{lemma:3} the kernel of $\1+K$ must be finite dimensional, $\Ker \Lambda_\pi R_\Pi R_{Q_\lambda}  \Lambda^{-1}_{\tilde q_\lambda}<\infty$; consequently 
\bq
\label{eq:56a}
r=\dim \Ker (R_\Pi R_{Q_\lambda}:\H(g,\tilde q_\lambda) \rightarrow H(g,\pi))< \infty.
\eq
 Next, let $U\subset \CT({\bf{X}}) \cap \H_\chi$ be a subspace such that
 \bqn
 A_\lambda[u] < 0, \qquad \forall \,\, 0\not=u \in U.
 \eqn
Then, for all $0\not= u \in V:= U\cap L\cap \complement_{\H(g,\tilde q _\lambda)}\,  (\Ker R_\Pi \,R_{Q_\lambda}:\H(g,\tilde q_\lambda) \rightarrow \H(g,\pi)) $,
 \bq
 \label{eq:56}
 0 > B_\lambda [R_\Pi \,R_{Q_\lambda}\, u].
 \eq
 Because  $R_\Pi\, R_{Q_\lambda}$ is injective on $V$,   \eqref{eq:56a} yields the inequality  $\dim U \leq \dim V+C \leq \dim R_\Pi\, R_{Q_\lambda} V +C$ for some constant $C>0$ independent of $\lambda$. 
Since $R_\Pi R_{Q_\lambda}$ commutes with the operators $T(k)$ of the representation of $G$,  $R_\Pi R_{Q_\lambda}V \subset \H_\chi \cap \H(g,\pi)$, and we obtain the estimate
 \bqn
 \dim U \leq \sup_{W \in \H(g,\pi) \cap \H_\chi} \mklm{\dim W: B_\lambda[w]< 0 \quad \forall \,0\not= w \in W } +C.
 \eqn
But $\CT(\R^n) \cap \H_\chi$ is dense in $\H(g,\pi) \cap \H_\chi$, and the assertion of the proposition follows.
 \end{proof}
 Let us now  prove Lemma \ref{lemma:10}. 
 
 \begin{proof}[Proof of Lemma \ref{lemma:10}]
 Let $ u \in \CT({\bf{X}})$. Then
 \bqn
 \Op^r(\eta_{c}) \, u (x) = \int \int e^{i(x-y)\xi} \eta_{c}(y,\xi) u (y) dy \, \dbar \xi =u(x),
 \eqn
 since $\eta_{c}$ is equal one on ${\bf{X}} \times \R^n$.
 Now, for general $B \in \Lcal (g,p)$, $\sigma^r( \Op^r(\eta_{c}) \, B)$ is given by an asymptotic expansion $\sum_j a_j$, where the first term is equal to $\eta_c \sigma^r(B)$. Consequently, $\sigma^r( \Op^r(\eta_c) B) =\eta_c \sigma^r(B) +(a- \eta_c \sigma^r(B)) +r$, with  $a$ as in Propostion \ref{prop:A}, and $r \in S^{-\infty}(h^{-2\delta}g,p)$. But  $a- \eta_c \sigma^r(B)=0$ on ${\bf{X}} \times \R^n$, and we obtain 
  \bq
  \label{eq:57}
  \Op^r(\eta_{c}) \, B \, u= B\,u +T\, u, \qquad T \in \Lcal^{-\infty}(h^{-2\delta}g,p).
 \eq 
Using Lemma \ref{lemma:11}, and setting $\tilde u =R_\Pi\, R_{Q_\lambda} u$, one computes
 \begin{align*}
 ((A_0-\lambda \1) \, u, \, u) & = ( \Op^w(\tilde a_\lambda) \,Q_\lambda \, \Pi \, \tilde u, \, Q_\lambda \, \Pi \, \tilde u) +(T_1 \, u, \, u)\\
 &=(\Pi^\ast [ Q_\lambda^\ast \,\Op^w(\tilde a_\lambda) \, Q_\lambda -4  \Op^w (h^{\delta-\omega}+ C_0 d)] \,\Pi \,\tilde u, \, \tilde u) \\
 &+4 ( \Pi^\ast  \Op^w(h^{\delta-\omega} +C_0 d) \Pi \,\tilde u , \, \tilde u ) + ( T_1 \,u,\, u)\\
 &=: ( \Pi^\ast \, D_1 \, \Pi \, \Op^r(\eta_{1}) \, \tilde u, \,\Op^r(\eta_{1}) \tilde u) + 4( D_2 \tilde u, \, \tilde u ) + (T_2 \, u, \, u),
 \end{align*}
 where we took \eqref{eq:57} into account together with $R_\Pi R_{Q_\lambda} -R_\Pi R_{Q_\lambda}   \in \Lcal^{-\infty}(g,\tilde q _\lambda \pi^{-1})$, and $T_i \in \Lcal^{-\infty}$. The reason for including $\Op^r(\eta_1)$ will become apparent in the proof of the next theorem. Now, by \eqref{eq:41a}, $a_\lambda-\sigma^w(Q_\lambda^\ast \Op^w (\tilde a_\lambda) Q_\lambda) \in S(g,d)$, which implies that for sufficiently large $C_0$ 
 \bqn
 D_1-\Op^w(a_\lambda^+)=Q_\lambda^\ast \Op^w(\tilde a_\lambda) Q_\lambda+ 4 C_0 \Op^w(d) -\Op^w(a_\lambda) \in \Lcal \I^+(g,d),
 \eqn
so that according to Lemma \ref{lemma:SI2}, there exists a $T_3 \in \Lcal^{-\infty}(g,d)$ such that $D_1 -\Op^w(a_\lambda^+) \geq T_3$. On the other hand, since $\pi^2  = (h^{\delta-\omega}+ C_0 d)^{-1}$, $D_2-\1 \in \Lcal(g,h)$, and we obtain
  \bq
  \label{eq:61ais}
  ( (A_0-\lambda \1) \, u, \, u ) \geq( \Op^l(\eta_{1})\, \Pi^\ast \, \Op^w(a_\lambda^+) \, \Pi \, \Op^r ( \eta_{1}) \, \tilde u , \, \tilde u ) + 2 \norm{\tilde u }^2 + (T_4 \, u , \, u ),
  \eq
  where $T_4 \in \Lcal (g, \pi^{-2} \tilde q_\lambda^2  h)$; hereby we used the fact that  $\Op^l(\eta_{1})$ is the adjoint of $\Op^r(\eta_{1})$, compare \cite{shubin}, page 26. Furthermore, since by \eqref{eq:12} the Weyl symbol of $R_\Pi R_{Q_\lambda}$ is equal to $\pi^{-1} \tilde q_\lambda$ modulo terms of lower order,
 \bqn
 (R_\Pi\, R_{Q_\lambda})^\ast R_\Pi\, R_{Q_\lambda} +T_4 \in \Lcal \I^+ (g , \pi^{-2} \tilde q_\lambda^2 ).
 \eqn
Lemmas \ref{lemma:SI1} -  \ref{lemma:3} now  allow us to deduce the  existence of a subspace $L \subset \CT({\bf{X}})$ of finite codimension in $\L^2({\bf{X}})$ such that 
 \bq
 \label{eq:61bis}
 \norm{\tilde u }^2 + (T_4 \, u, \, u)= ( [ (R_\Pi\, R_{Q_\lambda})^\ast R_\Pi\, R_{Q_\lambda} +T_4 ] u, \, u) >0
 \eq
 for all $0\not=u \in L$, and all $\lambda$. Indeed, according to Lemma \ref{lemma:SI2}, $ \Lambda_{\pi^2} [ (R_\Pi\, R_{Q_\lambda})^\ast R_\Pi\, R_{Q_\lambda} +T_4]\Lambda^{-1} _{\tilde q_\lambda^2} \in \mathcal{LI}^+(g,1)$  can be written in the form $B^\ast B +T_5$, where $B \in \mathcal{LI}(g,1)$ and $T \in \Lcal^{-\infty}(g,1)$. By a reasoning similar to the one that led to \eqref{eq:56a},  one can infer from Lemma \ref{lemma:SI1} that the kernel of $B$ must be finite dimensional, and together with Lemma \ref{lemma:3} conclude that there exists a subspace $\tilde L\subset \L^2(\R^n)$  of finite codimension such that 
 \bqn
  \norm{Bu}_{\L^2} \geq c \norm{u}_{\L^2}, \qquad  \norm{T_5u}_{L^2} < \frac {c^2} 2 \norm{u}_{\L^2},
 \eqn
 for all $u \in \tilde L$ and some constant $c>0$. Thus, we obtain  \eqref{eq:61bis}, and together with \eqref{eq:61ais} we get 
 \bqn
 ((A_0-\lambda \1) \, u, \,u) \geq ([\Op^l(\eta_{1})\, \Pi^\ast \Op^w(a_\lambda^+) \,\Pi \, \Op^r(\eta_{1})+\1 ] \, \tilde u, \, \tilde u ) 
 \eqn
 for all $0\not =u \in L$. This concludes the proof of the lemma. 
 \end{proof}
 We are now in position to prove an estimate from above for $\mathcal{N}(A_0-\lambda \1, \H_\chi\cap \CT({\bf{X}}))$.
 \begin{theorem}
 \label{thm:5}
 Let $M_\chi^{\F_\lambda}$ be the number of eigenvalues of $\F_\lambda$  which are $\geq 1/2$, and whose eigenfunctions are contained in the $\chi$-isotypic component $\H_\chi$. Then there exists a constant $C>0$ independent of $\lambda$ such that 
 \bqn 
 \mathcal{N}(A_0-\lambda \1, \H_\chi\cap \CT({\bf{X}})) \leq M_\chi^{\F_\lambda} +C.
 \eqn
 \end{theorem}
 \begin{proof}
We shall continue with the notation introduced in the proof of Proposition \ref{prop:8}. 
 According to that  proposition, it suffices to prove a similar estimate for $\Ncal(\Op^l(\eta_{1}) \Pi \Op^w(a_\lambda^+) \Pi \Op^r(\eta_{1}) +\1, \H_\chi\cap \CT(\R^n))$ from above. For this sake, we will show that there exists a subspace $L \subset U_\chi^{\F_\lambda}=\mathrm{Span} \mklm{ u \in \S(\R^n) \cap \H_\chi: \F_\lambda \, u =\nu \, u, \, \nu <1/2}$, whose codimension in $U_\chi^{\F_\lambda}$ is finite and uniformly bounded in $\lambda$, such that
 \bqn
 B_\lambda[u]\geq 0 \quad \text{for all } u \in L.
 \eqn
Indeed, let us assume this statement for a moment. Since $\F_\lambda$ is a compact self-adjoint operator in $\L^2(\R^n)$, there exists an orthonormal basis of eigenfunctions $\mklm{u_j}_{j=1}^\infty$ in $\S(\R^n)$. But $\F_\lambda $  commutes with the action $T(g)$ of $G$, so that each of the eigenspaces of $\F_\lambda$ is an invariant subspace, and must therefore decompose into a sum of irreducible $G$-modules. Consequently, 
  $\H_\chi$  has an orthonormal basis of eigenfunctions lying in $\S(\R^n) \cap \H_\chi$. Hence, 
 \bqn
 \H_\chi = U_\chi^{\F_\lambda} \oplus V_\chi^{\F_\lambda},
 \eqn
 where  $V_\chi^{\F_\lambda} = \mathrm{Span}\{ u \in \S(\R^n) \cap \H_\chi: \F_\lambda u=\nu u, \, \nu \geq 1/2\}$.
 Now, if  $W \subset \S(\R^n)\cap \H_\chi$ is a subspace with
 \bqn
 B_\lambda[u]<0 \quad \text{for all } 0 \not= u \in W,
 \eqn
 then $L \cap W =\mklm{0}$, and therefore $W \subset V_\chi^{\F_\lambda}\oplus U$, where $U$  is a finite dimensional subspace of $U_\chi^{\F_\lambda}$ whose dimension is bounded by some constant $C>0$ independent of $\lambda$.  Consequently, $\dim W \leq \dim V_\chi^{\F_\lambda} +C$. But this implies 
 \begin{align*}
 \Ncal((&\Op^l(\eta_{1}) \Pi \Op^w(a_{\lambda}^+) \Pi \Op^r(\eta_{1}) +\1), \H_\chi\cap \CT(\R^n))\\&\leq \sup _{W \subset \S(R^n) \cap \H_\chi} \mklm{ \dim W : ((\Op^l(\eta_{1}) \Pi \Op^w(a_{\lambda}^+) \Pi \Op^r(\eta_{1}) +\1) \, u , \, u)< 0 \quad \forall \,0\not = u \in W} \\
 &\leq \dim V_\chi^{\F_\lambda}+C =M_\chi^{\F_\lambda} +C,
 \end{align*}
 and the assertion of the theorem follows with the previous proposition. Let us now show the existence of the subspace $L$. Take $ v \in U_\chi^{\F_\lambda} \subset \L^2(\R^n)$, and put $\tilde v=(\1-\F_\lambda) v$. We then expect that $B_\lambda[\tilde v] \geq 0$. Now, one computes
 \begin{align}
 \label{eq:59}
 \begin{split}
 B_\lambda[\tilde v] &= ((\1-\F_\lambda') \Op^l(\eta_{1}) \Pi \Op^w(a_{\lambda}^+) \Pi \Op^r(\eta_{1})(\1-\F_\lambda') \, v, \, v) + \norm{(\1-\F_\lambda)v}^2 +(K_1 \, v, \, v)\\
 &\geq ( Dv,v) +(\norm{v}-\norm{\F_\lambda v})^2 + (K_1 v, v)\\
 & \geq (Dv, v) +\frac 1 4 \norm{v}^2 + (K_1 v,v), 
 \end{split}
 \end{align}
 where we put $\F_\lambda'= \Op^w(\chi_\lambda^+)^2(3-2\Op^w(\chi_\lambda^+))$, 
 \bqn
 D=(\1-\F_\lambda') \Op^l(\eta_{1}) \Pi \Op^w(a_{\lambda}^+) \Pi \Op^r(\eta_{1})(\1-\F_\lambda'),
 \eqn
 and  $K_1 \in \Lcal ^{-\infty}$.
  Indeed, one has $\norm{\F_\lambda v } \leq \frac 1 2 \norm{v}$, and $\Op^r(\eta_{1}) \F_\lambda - \Op^r(\eta_{1})\F_\lambda' \in \Lcal^{-\infty}(h^{-2\delta}g,1)$, 
 since the terms in the asymptotic expansions of the Weyl symbols of $\Op^r(\eta_{1}) \F_\lambda$ and $\Op^r(\eta_1)\F_\lambda'$  coincide because of $\eta_{2}=1$ on $\supp \eta_{1}$.  
Next we note that, similarly to \eqref{eq:44a},
\bq
\supp _\infty \sigma^w([\F_\lambda', \Op^l(\eta_{1}) \Pi]) \subset \supp_{\mathrm{diff}} \chi_\lambda^+ \subset \mklm{ (x,\xi): |a_\lambda^+ (x,\xi)| \leq h^\delta(x,\xi)},
\eq 
and we set
 \bqn
 b_\lambda^+= a_\lambda^+ \theta_\lambda, \qquad \theta_\lambda= \theta\Big (\frac 1 2 a_\lambda^+ h^{-\delta}\Big ), 
 \eqn
 with $\theta$ as in \eqref{eq:52c}. An argument similar to that concerning $b_\lambda^-$ shows that $b_\lambda^+ \in S(h^{-2\delta}g, h^\delta)$. Now, because of  $b_\lambda^+=a_\lambda^+$ on $\supp _\infty \sigma^w([\F_\lambda', \Op^l(\eta_{1} )\Pi ])$, we have
 \begin{align*}
 (D \, v,\, v) &= ( [(\1-\F_\lambda'), \Op^l(\eta_{1}) \Pi ] \Op^w (b_\lambda^+) \Pi \Op^r(\eta_{1})(\1- \F_\lambda') \, v, \, v)\\
 &+ (\Op^l(\eta_{1}) \Pi (\1 -\F_\lambda') \Op^w(a_\lambda^+) \Pi \Op^r(\eta_{1}) (\1-\F_\lambda') \, v, \, v) + (K_2 \, v, \, v),
  \end{align*}
 where $K_2$ is of order $-\infty$. Since $ [(\1-\F_\lambda'), \Op^l(\eta_{1}) \Pi ] \Op^w (b_\lambda^+) \Pi \Op^r(\eta_{1})(\1- \F_\lambda') \in \Lcal(h^{-2\delta} g, h^\delta \pi^2) \subset \Lcal(h^{-2\delta}g, h^\omega)$, we therefore obtain
 \bqn
  (D \, v, \, v) = (\Op^l(\eta_{1}) \Pi (\1-\F_\lambda') \Op^w(a_\lambda^+)  \Pi \Op^r( \eta_{1})(\1 -\F_\lambda') \, v, \, v) + (K_3 \, v , \, v),
  \eqn
 where $K_3 \in \Lcal(h^{-2\delta} g , h^\omega)$. Using a similar argument to commute $\Pi \Op^r(\eta_1)$ with $\1-\F_\lambda'$, we finally get
 \bq
 \label{eq:60}
 (D \, v, \, v) = (\Op^l(\eta_{1}) \Pi (\1-\F_\lambda') \Op^w(a_\lambda^+) (\1 -\F_\lambda') \Pi \Op^r( \eta_{1}) \, v, \, v) + (K_3 \, v , \, v),
 \eq
 where $K_3 \in \Lcal(h^{-2\delta} g , h^\omega)$. Now, the asymptotic expansion of the Weyl symbol of the operator $(\1-\F_\lambda') \Op^w(a_\lambda^+)(\1-\F_\lambda')$ gives
 \bq
 \label{eq:61}
 \sigma^w((\1-\F_\lambda') \Op^w(a_\lambda^+)(\1-\F_\lambda'))=[1-(\chi_\lambda^+)^2 (3-2\chi_\lambda^+)]^2 a_\lambda^+ + r
 \eq
with  $\supp_\infty r \subset \supp_{\mathrm{diff}} \chi_\lambda^+$. While computing $r$, we can therefore replace $a_\lambda^+$ by $b_\lambda^+$, so that $ r \in S(h^{-2\delta}g,h^\delta)$. As a consequence, \eqref{eq:60} and \eqref{eq:61}  yield
 \bqn
 (D \, v, \, v)= (\Op^l(\eta_{1}) \Pi \Op^w\big ([1-(\chi_\lambda^+)^2(3-2\chi_\lambda^+)]^2 a_\lambda^+ \big ) \Pi \Op^r(\eta_{1}) \, v, \, v) + ( K_4 \, v, \, v),
 \eqn
 where $K_4 \in \Lcal(h^{-2\delta}g, h^\omega)$. Hereby we used again the fact that  $\pi^2 h^\delta\sim h^\omega$. Next, one verifies that $[1-(\chi_\lambda^+)^2(3-2\chi_\lambda^+)]^2 a_\lambda^+ + C_1 h^\delta \in SI^+ (h^{-2\delta} g,[1-(\chi_\lambda^+)^2(3-2\chi_\lambda^+)]^2 a_\lambda^++C_1 h^\delta )$ for some $C_1>0$, since 
 $\chi_\lambda^+ =1$ for $a_\lambda ^+<0$, so that $[1-(\chi_\lambda^+)^2(3-2\chi_\lambda^+)]^2 a_\lambda^+\geq 0$. 
 According to Lemma \ref{lemma:SI2}, we therefore have
 \bqn
 \Op^w\big ([1-(\chi_\lambda^+)^2(3-2\chi_\lambda^+)]^2 a_\lambda^+ \big ) \geq K_5 \in \Lcal(h^{-2\delta}g, h^\delta),
 \eqn
 and we arrive at the estimate 
 \bqn
 (D \, v, \, v) \geq (K_6 \, v, \, v), \qquad K_6 \in \Lcal(h^{-2\delta} g, h^\omega).
 \eqn
  Together with \eqref{eq:59} we finally obtain the estimate
 \bqn
B_\lambda[\tilde v] \geq \frac 1 4 (v,v) + (K_7 \, v, \, v), \qquad K_7 \in \Lcal(h^{-2\delta} g, h^\omega).
 \eqn 
 Using the already familiar argument of Lemma \ref{lemma:3}, one infers the existence of a subspace $M \subset \L^2(\R^n)$ of finite codimension on which $ \1/4+K_7$ is positive definite. Putting $L:=(1-\F_\lambda) ( U_\chi^{\F_\lambda} \cap M) \subset U_\chi^{\F_\lambda}$ we therefore get 
 \bqn
 B_\lambda[w]  \geq 0, \qquad \text{for all } \, w \in L.
 \eqn
Furthermore, since $\1-\F_\lambda$ is injective on $U_\chi^{\F_\lambda}$, $\mathrm{codim}_{U_\chi^{\F_\lambda}}\,  L=\mathrm{codim}_{U_\chi^{\F_\lambda}}\,  (M\cap U_\chi^{\F_\lambda})\leq \mathrm{codim} \,M$, as desired. This completes the proof of the theorem.
  \end{proof}
    \begin{remark}
 The leading idea in  the proof of the last theorem was that each $v \in U^{\F_\lambda}_\chi$ has to be, approximately, an eigenvector of the corresponding spectral projection operator of $A$ with eigenvalue zero. For this reason, such a $v$ cannot satisfy $(Av,v) < \lambda \norm{v}^2$, nor be an element of $W$. 
 \end{remark}

\section{Asymptotics  for $\tr P_\chi\E_\lambda$ and $\tr P_\chi\F_\lambda$. The finite group case}

 For the rest of Part I, we shall concentrate on the case where $G$ is a finite group. The compact group case will be treated in Part II. The two preceding sections showed that, in view of  Lemmata \ref{lemma:12a} and \ref{lemma:12f}, the spectral counting function  $N_\chi(\lambda)=\mathcal{N}(A_0-\lambda \1,\H_\chi \cap \CT({\bf{X}}))$ can be estimated from below and from above in terms of the traces of  $P_\chi \E_\lambda$ and $P_\chi \F_\lambda$,  and their squares. We will therefore now proceed to estimate these traces in terms of the reduced Weyl volume. For this sake,  we introduce  first  certain sets associated to  the support of the symbols of the approximate spectral projection operators; their significance will become apparent later.  Thus, let 
 \begin{align*}
 \begin{split}
 W_\lambda&=\mklm{(x,\xi) \in {\bf{X}} \times \R^n: a_\lambda < 0}, \\
 A_{c,\lambda}&= \mklm{(x,\xi) \in {\bf{X}} \times \R^n: a_\lambda < c( h^{\delta-\omega}+d)},  \qquad B_{c,\lambda}= {\bf{X}} \times \R^n - A_{c,\lambda},\\
 D_c&=(\gd {\bf{X}}  \times \R^n) (c,h^{-2 \delta}g), \\
 F_\lambda&=\mklm{(x,\xi)\in {\bf{X}} \times \R^n: \chi_\lambda=0 \quad \text{or}  \quad \eta_{\lambda,-2} =0 \quad \text{or} \quad \chi_\lambda=\eta_{\lambda,-2} =1 },\\
 {\mathcal{RV}}_{c,\lambda}&=\mklm{(x,\xi) \in {\bf{X}} \times \R^n: | a_\lambda|< c(h^{\delta-\omega} + d)}\cup \mklm {(x,\xi) \in D_c: x \in {\bf{X}}, \, a_\lambda < c (h^{\delta -\omega} +d)}.
  \end{split}
  \end{align*}
 Note that $D_c=\mklm{ (x,\xi)\in \R^{2n}: \dist (x,\gd {\bf{X}}) <  \sqrt c\big (1 + |x|^2+|\xi|^2\big )^{-\delta/2}}$, since for 
 \bdm
 h^{-2\delta}(x,\xi)g_{(x,\xi)} ( x-y, \xi - \eta)= (1+|x|^2 + |\xi|^2)^\delta \Big [ \frac { |\xi-\eta|^2}{1 +  | x|^2 +|\xi|^2} + |  x -y|^2 \Big ] < c
 \edm
to hold  for some $(y,\eta) \in \gd {\bf{X}} \times \R^n $,  it is  necessary and sufficient  that  $|x-y|^2(1 +|x|^2+|\xi|^2)^\delta<c$ is satisfied for some $y \in \gd {\bf{X}} $. Now,  recall that $|G|= \sum_{\chi \in \hat G} d_\chi^2$. We then have the following 
  \begin{proposition}
  \label{prop:6a}
For sufficiently large $c>0$ we have
\bq
\label{eq:21}
|\tr P_\chi \E_\lambda - V_\chi({\bf{X}} \times \R^n,a_\lambda)| \leq c\, \vol \mathcal{RV}_{c,\lambda},
\eq
where
\bq
\label{eq:10}
 V_\chi({\bf{X}} \times \R^n,a_\lambda) =\frac {d_\chi^2}{|G|} \int \int_{{\bf{X}} \times  \R^n} \iota_{(-\infty, 0]}(a_\lambda(x,\xi))  dx \, \dbar \xi=\frac {d_\chi^2}{(2\pi)^n|G|} \,\vol W_\lambda 
\eq
is the expected approximation given in terms of the reduced Weyl volume, and $\iota_{(-\infty, 0]}$ denotes the characteristic function of the interval ${(-\infty, 0]}$. Furthermore, a similar estimate holds for $\tr \, P_\chi \E_\lambda \cdot P_\chi \E_\lambda$, too.
\end{proposition}
 \begin{proof}
 The proof will  require several steps. Let $\sigma^r(\E_\lambda)(x,\xi)$ denote the right symbol of $\E_\lambda$. Then, for $u \in \CT(\R^n)$, 
 \bqn
 P_\chi \E_\lambda u(x)= \frac {d_\chi}{|G|} \sum_{h \in G} \overline{\chi(h)}\int \int e^{i(h^{-1}  x -y ) \xi} \sigma^r(\E_\lambda) ( y, \xi ) u(y)  dy \, \dbar \xi.
 \eqn
 The kernel of $P_\chi\E_\lambda$, which is a rapidly decreasing function, is given by
 \bqn
 K_{P_\chi \E_\lambda}(x,y)= \frac {d_\chi}{|G|}\sum_{h \in G}   \overline{\chi(h)}\int e^{i(h^{-1}  x -y ) \xi} \sigma^r(\E_\lambda) ( y, \xi)    \, \dbar \xi.
 \eqn
The trace of $P_\chi \E_\lambda$ can therefore be computed by
 \begin{gather*}
 \tr P_\chi \E_\lambda=\int K_{P_\chi \E_\lambda}(x,x) dx \\
 =\frac {d_\chi^2}{|G|}  \tr \E_\lambda + \frac {d_\chi}{|G|} \sum _{h\not= e} \overline{\chi(h)} \int \int e^{i(h^{-1}x- x) \xi} \sigma^r(\E_\lambda) (x, \xi ) dx \, \dbar \xi,
 \end{gather*}
 where we made use of the relation $\chi(e)=d_\chi$, and the fact that $\tr \E_\lambda=\int \int \sigma^r(\E_\lambda) (x,\xi) dx \, \dbar \xi$. As a next step, we will prove that,  for all $e \not=h\in G$, there exists a sufficiently large constant $c>0$ such that
 \bq
 \label{eq:23a}
 \left | \int \int e^{i(h^{-1}x- x) \xi} \sigma^r(\E_\lambda) (x, \xi) dx \, \dbar \xi \right | \leq c \, \vol ( {\mathcal{RV}}_{c,\lambda}).
 \eq
 As already noticed, the decay properties of $\sigma^\tau(\E_\lambda)(x,\xi) \in S(h^{-2\delta}g,1)$ are independent of $\lambda$ for arbitrary $\tau \in \R$, while its support does depend on $\lambda$.  Indeed, by Theorem \ref{thm:MT} and Corollary \ref{cor:1}, together with the asymptotic expansions \eqref{eq:A} and \eqref{eq:12bis} and Proposition \ref{prop:A},  
\bq
\label{eq:23b}
\sigma^\tau(\E_\lambda) = (\eta_{\lambda,-2}^2 \chi_\lambda)^2 ( 3- 2 \eta_{\lambda,-2}^2\chi_\lambda) + f_\lambda + r_\lambda,
\eq
where $r_\lambda \in S^{-\infty}(h^{-2\delta}g,1)$, and $ f_\lambda \in S(h^{-2\delta} g, h^{1 - 2\delta})$,  everything uniformly in $\lambda$; in addition, $f_\lambda(x,\xi)=0$ if $(x,\xi) \in F_\lambda$. To see this, note that $\sigma^\tau(\E_\lambda)(x,\xi)$ is given asymptotically as a linear combination of products of derivatives of $\sigma^w( \E_\lambda)$ at $(x,\xi)$, which in turn is given asymptotically by a linear combination of terms involving derivatives of $\eta_{\lambda,-2},\, \chi_\lambda$. The $\tau$-symbol of  $\E_\lambda$ is therefore asymptotically given by 
\bqn
\sigma^\tau  ( \E_\lambda)  -\sum_{0 \leq j <N} a_j \in S(h^{-2\delta} g, h^{(1-2\delta)N}), \qquad a_j \in S(h^{-2\delta}g,h^{(1-2\delta)j}),
\eqn
where the first summand $a_0$ is equal to  $(\eta_{\lambda,-2}^2 \chi_\lambda )^2( 3- 2 \eta_{\lambda,-2}^2\chi_\lambda)$. Let now $a$ be as in Proposition \ref{prop:A} such that  $a\sim\sum_{j\geq 0} a_j$, and put  $r_\lambda=\sigma^\tau(\E_\lambda)-a \in S^{-\infty}(h^{-2\delta} g, 1)$. Since $\supp  \,(a -a_0) \subset \bigcup_{j\geq 1} \supp a_j$, $f_\lambda =a- (\eta_{\lambda,-2}^2 \chi_\lambda )^2( 3- 2 \eta_{\lambda,-2}^2\chi_\lambda)\in S(h^{-2\delta} g, h^{1 - 2\delta})$ must vanish on $F_\lambda$, and we obtain \eqref{eq:23b}. Now, since $|r_\lambda(x,\xi)| \leq C'(1+|x|^2+|\xi|^2)^{-N/2}$ for some constant $C'$ independent of $\lambda$ and $N$ arbitrarily large, we get the uniform bound 
\bdm
\int \int |r_\lambda(x,\xi)| dx \, \dbar \xi \leq C;
\edm
note that the $x$-dependence of $h(x,\xi)$ is crucial at this point. For this reason, and in order to show \eqref{eq:23a}, where now $\tau=1$, we can  restrict ourselves to the study of
 \bq
\label{eq:24}
 \int \int_{{\bf{X}} \times \R^n}  e^{i(h^{-1}x -x) \xi} ((\eta_{\lambda,-2}^2 \chi_\lambda )^2( 3- 2 \eta_{\lambda,-2}^2\chi_\lambda)+f_\lambda) (x ,\xi) dx \, \dbar \xi,
 \eq
 where we took into account that  $\eta_{\lambda,-2}$ has compact $x$-support in ${\bf{X}}$.
Next, we examine the geometry of the action of $G$ in more detail.  Thus, let
\bdm
\Sigma=\mklm{ x \in \R^n: gx =x \text{ for some } e\not=g \in G}
\edm
denote the set of not necessarily simultaneous fixed points of $G$. In other words,
 \bdm
 \Sigma = \bigcup _{e \not=g \in G} \Sigma_g, \qquad \Sigma_g = \mklm{x \in \R^n: gx =x}.
 \edm
 Note that
every connected component of $\Sigma_g$ is a closed, totally geodesic submanifold. We then have the following
\begin{lemma}
\label{prop:0}  
There exists a constant $\kappa>0$ such that $d(gx,x) \geq \kappa \, d(x,\Sigma_g)$ for all $x \in \R^n$,  and  arbitrary $e\not=g \in G$.
\end{lemma} 
\begin{proof}[Proof of Lemma \ref{prop:0}]
Let $x \in \R^n -\Sigma_g$ be an arbitrary point, and $p$ the closest point to $x$ belonging to $\Sigma_g$. Write $x = \exp_p t_0 X$, where $\exp_p$ denotes the exponential mapping of $\R^n$, and $(p,X) \in T_p(\R^n), \, |X|=1$. Then $t_0=\d(x,\Sigma_g)$.
Consider next the direct sum decomposition $T_p(\R^n)=U\oplus V$, where
\bdm
U= \mklm{(p,Y) \in T_p(\R^n): dg_p(Y)=Y},
\edm
and $V=U^\perp$. Since $p$ is a fixed point of $g$, we also have the identity
\bqn
g \, \exp_p Y= \exp_p dg_p (Y),
\eqn
which implies $\exp_p t Y \in \Sigma_g$ if, and only if, $(p,Y) \in U$, where $t \in \R$. Consequently, $U=T_p(\Sigma_g)$. Now, with $\exp _p t Y =p+tY$, and $x$, $p$ as above, 
one computes
\bdm
|gx -x|^2=|\exp_p t _0 dg_p(X) - \exp_p t_0 X|^2= |p + t_0 dg_p (X) - p -t_0X|^2=d^2(x,\Sigma_g) |dg_p(X)-X|^2.
\edm 
Because of $(x-p) \perp \Sigma_g$, we must have $(p,X)  \in T_p(\Sigma_g)^\perp=V$, and therefore $|dg_p(X) -X|^2 \not=0$.  The latter expression depends continuously on $(p,X) \in \mklm{(p,Y) \in T_p(\Sigma_g)^\perp: |Y|=1}$, and is actually independent of $p$, so that it can be estimated from below by some positive constant uniformly for all $x$. The assertion of the proposition now follows.
\end{proof}
Returning now to our previous computations, we split the integral in \eqref{eq:24} in an integral over
 \bdm
 \D=\mklm{(x,\xi) \in {\bf{X}} \times \R^n: \dist (x, \Sigma) \geq  ( 1+ |\xi|^2)^{-\delta/2}},
  \edm
  and a second integral over the complement of $\D$ in ${\bf{X}} \times \R^n$. Since $\supp \chi_\lambda \subset \{ (x,\xi): a_\lambda +4 h^{\delta-\omega} + 8 C_0 \leq h^\delta \}$, the integral over $\complement_{\Omega \times {\bf{X}}} \D$ can be estimated by a constant times the volume of the set $\{(x,\xi) \in \complement_{\Omega \times \bf{X}} \D: a_\lambda +4h^{\delta-\omega}+8 C_0 d \leq h^\delta\}$, which is contained in the set $\{(x,\xi) \in {\bf{X}} \times \R^n: \dist (x,\Sigma) <  (1+|\xi|^2)^{-\delta/2}, a_\lambda \leq  c(h^{\delta-\omega}+d)\}$ for some sufficiently large $c>0$. By examining the proof of Lemma \ref{lemma:13}, one sees that the volume of the latter can be estimated from above by
  \bqn
  \int_{K \leq |\xi| \leq c_1 \lambda^{1/2m}} \vol (\Sigma_{c_2 |\xi|^{-\delta}} \cap {\bf{X}}) \, d\xi+c_3
  \eqn
 for some suitable constants $K,c_i>0$, and consequently has the same asymptotic behaviour in $\lambda $ as the volume of $\mathcal{RV}_{c,\lambda}$. In studying the asymptotic behaviour of the integral \eqref{eq:24}, we can therefore restrict the domain of integration to $\D$. By the previous lemma, there exists a constant $\kappa>0$ such  that 
 \bqn
 |h^{-1}x-x| \geq  \kappa ( 1+|\xi|^2)^{-\delta/2} \quad \text{for all } (x,\xi)  \in \D \text{ and } e\not=h \in G.
 \eqn
 Since $(\eta_{\lambda,-2}^2 \chi_\lambda )^2( 3- 2 \eta_{\lambda,-2}^2\chi_\lambda)+f_\lambda$ has compact support in $\xi$, this implies that 
 \bdm
 \frac{e^{i(h^{-1}x-x) \xi}}{|h^{-1}x-x|^2} \gd_\xi^\alpha ((\eta_{\lambda,-2}^2 \chi_\lambda )^2( 3- 2 \eta_{\lambda,-2}^2\chi_\lambda)+f_\lambda)(x, \xi)
 \edm
 is integrable on $\D$, as well as rapidly decreasing in $\xi$. Integrating  by parts with respect to $\xi$ we  therefore get for  \eqref{eq:24} the expression
 \bq
 \label{eq:26}
 \int \int _\D \frac{e^{i(h^{-1}x-x) \xi}}{|h^{-1}x-x|^2} (-\gd^2_{\xi_1}- \cdots -\gd^2_{\xi_n} ) ((\eta_{\lambda,-2}^2 \chi_\lambda )^2( 3- 2 \eta_{\lambda,-2}^2\chi_\lambda)+f_\lambda)(x, \xi)d x\,  \, \dbar \xi;
 \eq
  in particular notice that, by Fubini's Theorem,  the boundary contributions vanish. Now, if $(x,\xi)\in  F_\lambda$, the function $(\eta_{\lambda,-2}^2 \chi_\lambda )^2( 3- 2 \eta_{\lambda,-2}^2\chi_\lambda)+f_\lambda$ is constant, so its derivatives with respect to $\xi$ are zero, and we can restrict the integration in \eqref{eq:26} to the set $\complement_{{\bf{X}} \times \R^n}  F_\lambda\cap \D$, where $\complement_{{\bf{X}} \times \R^n}  F_\lambda$ denotes the complement of $\F_\lambda$ in ${\bf{X}} \times \R^n$.
 \begin{lemma}
 For sufficiently large $c>0$, the set $\complement_{{\bf{X}} \times \R^n} F_\lambda$ is contained in ${\mathcal{RV}}_{c,\lambda}$.
\end{lemma}
\begin{proof}
This assertion is already stated in \cite{levendorskii}, page 55. For the sake of completeness, we give a proof here.  Thus, consider 
 $$E_\lambda= \mklm{ (x,\xi) \in {\bf{X}} \times \R^n:  (x,\xi) \not\in D_4, a_\lambda < -4h^{\delta-\omega} - 8 C_0 d },$$ 
 and let $\M_\lambda$ be defined as in \eqref{eq:19}. 
 Since $\supp \tilde \eta_2 \subset D_4$, and $\psi_{\lambda,1/2} =1$ on $\M_\lambda(1/2,h^{-2\delta}g)$, it is clear that  
\begin{align}
\label{eq:27}
\begin{split}
E_\lambda \subset \mklm{(x,\xi) \in {\bf{X}}\times \R^n: \chi_\lambda=\eta_{\lambda,-2}=1} \subset F_\lambda,
\end{split}
\end{align}
and consequently $\complement_ {{\bf{X}} \times \R^n} F_\lambda\subset \complement_{{\bf{X}} \times \R^n}E_\lambda$. Next, we are going to prove that, for sufficiently large $c$, $(x,\xi) \in B_{c,\lambda}$ implies $(x,\xi) \not \in \M_\lambda(1, h^{-2\delta}g)$. Thus, assume  $(x,\xi) \in B_{c,\lambda}$; on ${\bf{X}}_\epsilon \times \mklm{\xi: |\xi| >1}$ we have
\bdm
c \Big ( \frac 1 {|\xi|} + \frac 1 {( 1 +|x|^2 + |\xi|^2 ) ^{(\delta-\omega)/2)}}\Big ) \leq \frac {|\xi|^{2m}} {|\xi|^{2m}+\lambda} \Big (1 - \frac \lambda {a_{2m}(x,\xi)}\Big ). 
\edm
Therefore, as $c$ becomes large, $|\xi|$ must become large, too. On the other hand, if $(y,\eta)\in \M_\lambda$, $|\eta|$ must be bounded. For large $c$ we therefore have $|\xi-\eta|^2 \sim  |\xi|^2$, which means that $h^{-2\delta}(x,\xi) g_{(x,\xi)}(x-y,\xi-\eta)\sim(1+ |x|^2+|\xi|^2)^\delta \to \infty$ as $c\to \infty$. Hence, for sufficiently large $c$, $(x,\xi) \not\in \M_\lambda(1,h^{-2\delta}g)$. Since $\supp \psi_{\lambda,1/2} \subset \M_\lambda(1,h^{-2\delta}g)$, we arrive in this case at the inclusions 
\bq
\label{eq:28}
B_{c,\lambda} \subset \mklm{(x,\xi) \in {\bf{X}} \times \R^n: \eta_{\lambda,-2}(x,\xi)=0} \subset F_\lambda,
\eq
and combining \eqref{eq:27} and \eqref{eq:28} we get 
\bq
\label{eq:29}
\complement_{{\bf{X}} \times \R^n} F_\lambda \subset A_{c,\lambda} \cap \complement _{{\bf{X}}\times \R^n} E_\lambda \subset {\mathcal{RV}}_{c,\lambda},
\eq
as desired.
\end{proof}
As a consequence of the foregoing lemma,  the integral in \eqref{eq:26} is bounded from above by  the volume of ${\mathcal{RV}}_{c,\lambda}$, times a constant independent of $\lambda$, since the integrand is uniformly bounded with respect to  $\lambda$. Thus, we have shown  \eqref{eq:23a}. The assertion of the Proposition now follows by observing that
\bq
\label{eq:30}
\Big | \tr \E_\lambda - \frac {\vol W_\lambda}{(2\pi)^n} \Big | \leq c \, \vol\mathcal{RV}_{c,\lambda}.
\eq
Indeed,  similarly to  our previous discussion of the integral $\int \int e^{i(h^{-1}x-x)\xi} \sigma^r(\E_\lambda)(x,\xi) dx \, \dbar \xi$, the integral
 $$\tr \E_\lambda=\int \int \sigma^r(\E_\lambda)(x,\xi) dx \, \dbar \xi=\int \int ((\eta_{\lambda,-2}^2 \chi_\lambda )^2( 3- 2 \eta_{\lambda,-2}^2\chi_\lambda) + f_\lambda + r_\lambda)(x,\xi) dx \, \dbar \xi$$
  can be split into three parts; the contribution coming from $r_\lambda(x,\xi)$
is bounded by some constant independent of $\lambda$, while the contribution coming from $f_\lambda$ can be estimated in terms of the volume of $\mathcal{RV}_{c,\lambda}$, since $\supp f_\lambda \subset \complement_{{\bf{X}} \times \R^n} F_\lambda \subset {\mathcal{RV}}_{c,\lambda}$, by the previous lemma. Now, 
 $(\eta_{\lambda,-2}^2 \chi_\lambda )^2( 3- 2 \eta_{\lambda,-2}^2\chi_\lambda)$ must be equal  $1$ on $W_\lambda \cap \complement _{{\bf{X}}\times \R^{n}}  {\mathcal{RV}}_{c,\lambda}$, since according to \eqref{eq:29} we have $\complement_{{\bf{X}} \times \R^n}{\mathcal{RV}}_{c,\lambda} \subset B_{c,\lambda} \cup E_\lambda$, and hence $W_\lambda \cap \complement _{{\bf{X}}\times \R^{n}}  {\mathcal{RV}}_{c,\lambda} \subset E_\lambda \subset \mklm {(x,\xi) \in {\bf{X}} \times \R^n: \chi_\lambda =\eta_{\lambda,-2} =1}$, due to the fact that $W_\lambda \cap B_{c,\lambda}=\emptyset$. Furthermore, $(\eta_{\lambda,-2}^2 \chi_\lambda)^2 ( 3- 2 \eta_{\lambda,-2}^2 \chi_\lambda)$ vanishes on $B_{c,\lambda}$, since for large $c$, $(x,\xi) \in B_{c,\lambda}$ implies $(x,\xi) \not \in \M_\lambda(1, h^{-2\delta}g)$, by the proof of the previous lemma.  Taking into account that $ W_\lambda$ and ${\mathcal{RV}}_{c,\lambda}$ are subsets of $A_{c,\lambda}$,  we therefore obtain for sufficiently large $c$
 \begin{gather*}
 \int \int ((\eta_{\lambda,-2}^2 \chi_\lambda)^2 ( 3- 2 \eta_{\lambda,-2}^2 \chi_\lambda))(x,\xi) dx \, \dbar \xi\\ =\frac{\vol \big (W_\lambda \cap \complement _{A_{c,\lambda}}  {\mathcal{RV}}_{c,\lambda}\big )}{(2\pi)^n}+ \int \int_{A_{c,\lambda} - (W_\lambda \cap \complement _{A_{c,\lambda}}  {\mathcal{RV}}_{c,\lambda})}  ((\eta_{\lambda,-2}^2 \chi_\lambda)^2 ( 3- 2 \eta_{\lambda,-2}^2 \chi_\lambda))(x,\xi) dx \, \dbar \xi.
 \end{gather*}
Now, since $\complement _{A_{c,\lambda}} {\mathcal{RV}}_{c,\lambda} \subset W_\lambda$, one has  $A_{c,\lambda} -W_\lambda \cap \complement _{A_{c,\lambda}}  {\mathcal{RV}}_{c,\lambda}= {\mathcal{RV}}_{c,\lambda} $.
 The estimate \eqref{eq:30} now follows, and together with  \eqref{eq:23a} we obtain \eqref{eq:21}. Finally, if  in the previous computations  $\E_\lambda$ is replaced by $\E_\lambda^2$, we obtain a similar estimate for the trace of $P_\chi \E_\lambda\cdot P_\chi \E_\lambda=P_\chi \E_\lambda^2$. This concludes the proof of the proposition.
 \end{proof}
 As a consequence, we get the following
\begin{theorem}
\label{thm:2}
Let $N^{\E_\lambda}_\chi$ be the number of eigenvalues of $\E_\lambda$ which are $\geq 1/2$ and whose eigenfunctions are contained in the $\chi$-isotypic component $\H_\chi$ of $\L^2(\R^n)$. Then
\bq
\label{eq:38}
|N^{\E_\lambda}_\chi-V_\chi({\bf{X}} \times \R^n,a_\lambda)| \leq  c\, \vol\mathcal{RV}_{c,\lambda}
\eq
for some sufficiently large $c>0$.
\end{theorem}
 \begin{proof}
From the preceding proposition, and  the estimate \eqref{eq:12a}, one deduces that for some sufficiently large $c>0$ 
 \begin{align*}
 N^{\E_\lambda}_\chi&\leq 3 \tr P_\chi \E_\lambda  -2 \tr P_\chi \E_\lambda \cdot P_\chi \E_\lambda +c_2 \leq   V_\chi({\bf{X}} \times \R^n, a_\lambda) +c\, \vol \mathcal{RV}_{c,\lambda},\\
 N^{\E_\lambda}_\chi&\geq 2 \tr P_\chi \E_\lambda \cdot P_\chi \E_\lambda -\tr P_\chi \E_\lambda  -c_1\geq V_\chi({\bf{X}} \times \R^n, a_\lambda) -c \,\vol \mathcal{RV}_{c,\lambda},
 \end{align*}
 which completes the proof of \eqref{eq:38}.
 \end{proof}
In analogy to the previous considerations, one  proves the following
\begin{theorem}
\label{thm:2a}
For sufficiently large $c>0$ one has the estimate
\bqn
|M_\chi^{\F_\lambda} -V_\chi({\bf{X}} \times \R^n,a_\lambda)| \leq c\, \vol \mathcal{RV}_{c,\lambda},
\eqn
where $M_\chi^{\F_\lambda}$ is the number of eigenvalues of $\F_\lambda$, counting multiplicities, greater or equal $1/2$, and whose eigenfunctions are contained in the  $\chi$-isotypic component $\H_\chi$ of $\L^2(\R^n)$.
\end{theorem}
\begin{proof}
The proof   is similar to the one of Theorem \ref{thm:2}, and uses Lemma \ref{lemma:12f}.   In particular, 
 as in  equation \eqref{eq:23b}, one has
\bq
\label{eq:23z}
\sigma^\tau(\F_\lambda) = (\eta_{2}^2 \chi^+_\lambda)^2 ( 3- 2 \eta_{2}^2\chi^+_\lambda) + f_\lambda + r_\lambda,
\eq
where $r_\lambda \in S^{-\infty}(h^{-2\delta}g,1)$, and $ f_\lambda \in S(h^{-2\delta} g, h^{1 - 2\delta})$,  everything uniformly in $\lambda$. The aymptotic analysis for  $\tr P_\chi \F_\lambda$ and $\tr(P_\chi \F_\lambda)^2$ now follows the lines of the proof of Proposition \ref{prop:6a}.
\end{proof}

 \section{Proof of Theorem \ref{thm:1}}

We collect all the results obtained so far in the following  
 \begin{proposition}
\label{prop:6}
 There exist constants $C_1, C_2>0$ which do not depend on  $\lambda$, such that for all $\lambda$
 \bqn
 |  \mathcal{N}(A_0-\lambda \1, \H_\chi\cap \CT({\bf{X}})) - V_\chi({\bf{X}} \times \R^n, a_\lambda)| \leq C_1 \vol \mathcal{RV}_{C_1,\lambda}+C_2.
 \eqn
 \end{proposition} 
 \begin{proof}
By  Theorems \ref{thm:4} and \ref{thm:5}, there exist constants $C_i>0$ independent of $\lambda$ such that  $N^{\E_\lambda}_\chi -C_1 \leq  \mathcal{N}(A_0-\lambda \1, \H_\chi\cap \CT({\bf{X}})) \leq M_\chi^{\F_\lambda} +C_2$. Theorems \ref{thm:2} and \ref{thm:2a} then yield the estimate
 \bqn
 -c\,\vol \mathcal{RV}_{c,\lambda} - C_1 \leq  \mathcal{N}(A_0-\lambda \1, \H_\chi\cap \CT({\bf{X}})) - V_\chi({\bf{X}} \times \R^n, a_\lambda) \leq c \, \vol \mathcal{RV}_{c,\lambda} + C_2
 \eqn
 for some sufficiently large $c>0$.
  \end{proof}
 In order to formulate the main result, we need two last lemmata.
 \begin{lemma}
 \label{lemma:12}
Let $\gamma = \frac 1 n \int_{\bf{X}} \int _{S^{n-1}} \big ( a_{2m}(x,\xi) \big ) ^{-n/2m} dx \, \dbar \xi$, where $2m$ is the order of $A_0$. Then
\bqn
V_\chi({\bf{X}} \times \R^n, a_\lambda) = \frac {d_\chi^2}{|G|} \gamma \cdot \lambda^ {n/2m} +C
\eqn 
for some constant $C>0$ independent of $\lambda$.
 \end{lemma}
 \begin{proof}
 The reduced Weyl volume was defined in \eqref{eq:10} as 
 \bqn
 V_\chi({\bf{X}} \times \R^n, a_\lambda) = \frac {d_\chi^2}{|G|} \int \int _{{\bf{X}} \times \R^n} \iota_{(-\infty, 0]} (a_\lambda(x,\xi)) dx \,\dbar \xi =\frac {d_\chi^2}{(2\pi)^n|G|} \vol W_\lambda,
 \eqn
 where $W_\lambda=\mklm{ (x,\xi) \in {\bf{X}} \times \R^n: a_\lambda(x,\xi) < 0}$. Now, for some sufficiently small $\rho>0$, on ${\bf{X}}_\rho \times \mklm{\xi: |\xi| >1}$, $a_\lambda$ is given  by
 \bqn
 a_\lambda(x,\xi) = \frac 1 { 1  + \lambda|\xi|^{-2m}} \Big ( 1 - \frac \lambda {a_{2m}(x,\xi)}\Big ).
 \eqn
By \cite{levendorskii}, Lemma 13.1, condition \eqref{eq:8a} implies that $a_{2m}(x,\xi) \geq \iota >0$ for all $(x,\xi) \in {\bf{X}} \times S^{n-1}$, so that $a_\lambda$ is strictly negative on ${\bf{X}}_\rho \times \mklm{\xi: |\xi| >1}$ if, and only if, $a_{2m}(x,\xi)-\lambda <0$, which in turn is equivalent to 
 \bqn
 |\xi| < [ \lambda  \, a_{2m}^{-1}(x,\xi/|\xi|)]^{1/2m},
 \eqn
due to the homogeneity of the principal symbol. From this one concludes
\begin{align*}
V_\chi({\bf{X}} \times \R^n, a_\lambda)&=\frac{d_\chi^2}{(2\pi)^n|G|} \big [ \vol \mklm{(x,\xi) \in {\bf{X}} \times \R^n: |\xi|\leq 1}\\&+\vol \mklm{(x,\xi) \in {\bf{X}} \times \R^n: |\xi|^{2m} < \lambda \, a_{2m}^{-1}(x, \xi/|\xi|)}\big ]\\
&=O(1) + \frac{d_\chi^2}{(2\pi)^n|G|} \int_{\bf{X}} \int_{S^{n-1}} \int _0 ^{( \lambda/a_{2m}(x,\eta))^{1/2m}} r^{n-1} dr \, dS^{n-1}(\eta) dx\\
&= O(1) + \frac{d_\chi^2}{(2\pi)^n|G|} \int_{\bf{X}} \int _{S^{n-1}} \frac 1 n ( \lambda/a_{2m}(x,\eta))^{n/2m} dS^{n-1}(\eta) dx.
  \end{align*}
 \end{proof}
 \begin{lemma}
 \label{lemma:13} Assume that for some sufficiently small $\rho >0$ there exists a constant $C>0$ such that $\vol (\gd {\bf{X}})_\rho \leq C \rho$. Then 
  $\vol \mathcal{RV}_{c,\lambda}=\mathrm{O}(\lambda^{(n-\epsilon)/2m})$, where $\epsilon \in (0, \frac 1 2)$.
 \end{lemma}
 \begin{proof}
 According to the definition of $\mathcal{RV}_{c,\lambda}$ at the beginning of Section 6, we have 
\begin{align*}
\vol \mathcal{RV}_{c,\lambda}&\leq \vol \mklm{(x,\xi) \in {\bf{X}} \times \R^n:|a_\lambda| - c ( h^{\delta -\omega} + d)<0} \\
&+\vol \mklm{(x,\xi) \in D_c: x \in {\bf{X}}, a_\lambda < c(h^{\delta-\omega}+d)},
\end{align*}
where $D_c=\{ (x, \xi): \mathrm{dist}(x,\gd {\bf{X}}) < \sqrt c  ( 1 +  |x|^2 + |\xi|^2) ^{-\delta/2}\}$, and $0 < \delta -\omega < 1/2$. In what follows, let us assume that $\lambda \geq 1$.
It is not difficult to see that, for $|\xi|>1$, there exists a constant $c_1 >0$ independent of $\lambda$ such that 
\bq
\label{eq:70}
 a_\lambda(x,\xi)-c(h^{\delta-\omega}+ d)(x,\xi) < 0  \quad \Longrightarrow \quad |\xi| < c_1 \lambda^ {1/2m}.
\eq
Indeed, let $c_1$ be such that 
\bqn
c_1^{2m} \geq \max \big (2,  2/\iota \big ), \qquad \sup_{x \in {\bf{X}}, |\xi|>c_1} c(h^{\delta-\omega}+d) (x,\xi) \leq \frac 13,
\eqn
where $\iota>0$ is a lower bound for  $a_{2m}(x,\xi)$ on  ${\bf{X}} \times S^{n-1}$.
Since 
\bqn
1- \frac \lambda {a_{2m}(x,\xi)} \geq \frac 1 2 \quad \Longleftrightarrow \quad  |\xi |^{2m} \geq \frac {2 \lambda} {a_{2m}(x, \xi/|\xi|)},
\eqn
one computes for $|\xi| \geq c_1 \lambda^{1/2m}$ that
\bqn
a_\lambda(x,\xi)\geq  \frac 1 2 \frac 1{1+ \lambda |\xi|^{-2m}} \geq \frac 1 2 \frac 1 {1+c_1^{-2m}} \geq \frac 1 3,
\eqn
while, on the other hand, $c(h^{\delta-\omega} +d)(x,\xi) \leq  \frac 1 3$, so that $a_\lambda(x,\xi) -c(h^{\delta -\omega} +d)(x,\xi) \geq 0$. This proves \eqref{eq:70}. As a consequence, we obtain the estimate
\begin{align*}
\vol \{(x,\xi) &\in D_c:  x \in {\bf{X}}, a_\lambda < c(h^{\delta-\omega}+d)\} \\
&\leq \vol\mklm{ (x, \xi)\in {\bf{X}} \times \R^n:  |\xi| \geq K,  \mathrm{dist}(x, \gd {\bf{X}}) < c_2 |\xi|^{-\delta}, a_\lambda < c(h^{\delta -\omega} +d ) }\\
&+  \vol \mklm{(x,\xi) \in {\bf{X}} \times \R^n: |\xi| \leq K}\\
& \leq  \int_{K \leq |\xi| \leq c_1 \lambda ^{1/2m}} \vol \big ( (\gd {\bf{X}})_{c_2 |\xi|^{-\delta}}\cap {\bf{X}}\big ) \, d\xi+ c_3,
\end{align*}
where $\delta \in ( \frac 1 4, \frac 12)$, and $K\geq 1$ is some sufficiently large constant; here and in all what follows, $c_i>0$ denote suitable, positive constants independent of $\lambda$. Now, since $\vol (\gd {\bf{X}})_\rho \leq C \rho$, for some  $\rho >0$, 
\begin{align*}
\vol \{(x,\xi) &\in D_c:  x \in {\bf{X}}, a_\lambda < c(h^{\delta-\omega}+d)\} \leq c_2 \,  c_4 \int_{S^{n-1}}\int_{K \leq r \leq c_1 \lambda^ {1/2m}} r^{n-1-\delta} dr dS^{n-1}(\eta)  + c_3\\
&=c_5 (\lambda^{(n-\delta)/2m}-K^{n-\delta})+c_3.
\end{align*}
Next, let $|\xi| \geq K$, and assume that  the inequality
$|a_\lambda(x,\xi) | \leq c( h^{\delta -\omega}+d)(x,\xi)$ is fulfilled. As before, we have $|\xi|^{2m}<c_1^{2m} \lambda$, as well as
\bq
\label{eq:71}
\Big | 1 -\frac \lambda {a_{2m}(x,\xi)}\Big |  \leq c(1+\lambda|\xi|^{-2m} ) (d + h^{\delta-\omega}) (x,\xi)\leq c_6 (1+\lambda|\xi|^{-2m})|\xi|^{-(\delta -\omega)}. 
\eq
Combining \eqref{eq:70} and \eqref{eq:71}, one deduces for sufficiently large $K$ that  
 \bqn
 |\xi|^{2m} \geq -c_6 (|\xi|^{2m}+\lambda)|\xi|^{-(\delta -\omega)}+ \frac \lambda {a_{2m}(x,\xi/|\xi|)}\geq c_7 \lambda.
 \eqn
 Let us now introduce the variable $R(x,\xi)=  \lambda/a_{2m}(x,\xi)=\lambda|\xi|^{-2m}/a_{2m}(x,\xi/|\xi|)$.  Performing the corresponding change of variables one computes  
 \begin{align*}
 \vol &\{(x,\xi)  \in {\bf{X}} \times \R^n: |a_\lambda| -c ( h^{\delta-\omega}+d) <0\}\\
 &\leq \vol \mklm{ (x,\xi) \in {\bf{X}} \times \R^n:c_1 \lambda^{1/2m} >  |\xi|\geq K,\,  | 1 -R (x,\xi)|   \leq c_6 (1+\lambda|\xi|^{-2m})|\xi|^{-(\delta -\omega)} }+c_{8}\\
  & \leq \int_{\bf{X}} \int_{S^{n-1}} \int_{\{r \geq K: |1-R|\leq c_9 \lambda^{-(\delta -\omega)/2m}\}} r^{n-1} 
 \,dr \,dS^{n-1}(\eta) \, dx+ c_{8}\\
 &  \leq c_{10}  \int_{\bf{X}} \int_{S^{n-1}} \int_{\{R: |1-R|\leq c_9 \lambda^{-(\delta -\omega)/2m}\}} R^{-1} \left ( \frac \lambda{ R a_{2m}(x,\eta)} \right ) ^{\frac {n}{2m}}    \,dR\,dS^{n-1}(\eta) \, dx+ c_{8}\\
 &  \leq c_{11} \lambda^{\frac n{2m}} \int_{\{R: |1-R|\leq c_9 \lambda^{-(\delta -\omega)/2m}\}} R^{-\frac {n}{2m}-1}    \,dR+ c_{8}= O( \lambda^ {(n-(\delta -\omega))/2m}) +c_8.
 \end{align*}
 Hereby we made use of the fact that $(1+z)^\beta-(1-z)^\beta=O(|z|)$ for arbitrary $ z \in \C$, $|z|<1$, and $\beta \in \R$.
  \end{proof}
We are now in position to prove the main result of  Part I, which generalizes Theorem 13.1 of \cite{levendorskii} to bounded domains with symmetries in the finite group case.
\begin{theorem1}
Let $G$ be a finite group of isometries in Euclidean space $\R^n$, and ${\bf{X}} \subset \R^n$ a bounded domain which is invariant under $G$  such that, for some sufficiently small $\rho >0$, $\vol ( \gd {\bf{X}} )_\rho \leq C \rho$. Let further $A_0$ be a  symmetric, classical pseudodifferential operator in $\L^2(\R^n)$ of order $2m$ with $G$-invariant Weyl symbol $\sigma^w(A_0) \in S(g,h^{-2m})$ and principal symbol $a_{2m}$, and assume that $(A_0 \, u, \, u) \geq c \norm{u}^2_m$ for some $c>0$ and  all $u \in \CT({\bf{X}})$. Consider further  the Friedrichs extension of the operator
\bqn
\mathrm{res} \circ A_0 \circ \mathrm{ext}: \CT({\bf{X}}) \longrightarrow \L^2({\bf{X}}),
\eqn
and denote it by $A$. Finally, let $N_\chi(\lambda)$ be  the number of  eigenvalues of $A$ less or equal $\lambda$ and  with eigenfunctions in the  $\chi$-isotypic component $\mathrm{res}\, \H_\chi$ of $\L^2({\bf{X}})$, if $(-\infty,\lambda)$ contains no points of the essential spectrum, and equal to $\infty$, otherwise. Then, for all $\epsilon \in (0,\frac 1 2)$, 
\begin{equation*}
N_\chi(\lambda) =\frac {d^2_\chi}{|G|}\gamma   \lambda^{n/2m} +O(\lambda^{(n-\epsilon)/2m}),
\end{equation*}
where $d_\chi$ denotes the dimension of the  irreducible representation of $G$ corresponding to the character $\chi$, and 
\bdm
\gamma=\frac 1 n \int_{\bf{X}} \int_{S^{n-1}} (a_{2m}(x,\xi))^{-n/2m} dx \, \dbar \xi.  
\edm
In particular, $A$ has discrete spectrum.
\end{theorem1}
\begin{proof}
By Lemma \ref{lemma:8} and Proposition \ref{prop:6} we have
\bqn
|N_\chi(\lambda) -V_\chi({\bf{X}} \times \R^n, a_\lambda) | \leq C_1 \vol \mathcal{RV}_{C_1,\lambda}+ C_2
\eqn
for some suitable constants $C_1,C_2>0$ independent of $\lambda$. Lemma \ref{lemma:12} and \ref{lemma:13} then imply
\bqn
-\mathrm{O}(\lambda^{(n-\epsilon)/2m}) \leq N_\chi(\lambda) - \frac{d_\chi^2}{|G|} \gamma \lambda^{n/2m} \leq \mathrm{O}(\lambda^{(n-\epsilon)/2m})
\eqn
with  arbitrary $\epsilon \in (0,1/2)$. In particular, $N_\chi(\lambda)$ remains finite for $\lambda<\infty$, so that the essential spectrum of $A$ must be empty. The assertion of the theorem now follows.
\end{proof}

\providecommand{\bysame}{\leavevmode\hbox to3em{\hrulefill}\thinspace}
\providecommand{\MR}{\relax\ifhmode\unskip\space\fi MR }
\providecommand{\MRhref}[2]{%
  \href{http://www.ams.org/mathscinet-getitem?mr=#1}{#2}
}
\providecommand{\href}[2]{#2}



\end{document}